\documentclass[10pt]{amsart}
\usepackage{amscd,amssymb,graphics}

\usepackage{amsfonts}
\usepackage{amsmath}
\usepackage{amsxtra}
\usepackage{latexsym}
\usepackage[mathcal]{eucal}

\usepackage{graphics,colortbl}

\usepackage{amsmath,amssymb,amsthm}
\newtheorem*{example*}{Example}
\newtheorem*{remark*}{Remark}

\input xy
\xyoption{all}
\usepackage{epsfig}

\usepackage[pdftex,bookmarks,breaklinks]{hyperref}

\newtheorem{theorem}{Theorem}[section]
\newtheorem{fact}[theorem]{Fact}
\newtheorem{corollary}[theorem]{Corollary}
\newtheorem{lemma}[theorem]{Lemma}
\newtheorem{proposition}[theorem]{Proposition}

\newtheorem{question}[theorem]{Question}
\newtheorem{definition}[theorem]{Definition}
\numberwithin{equation}{section}
 
\newtheorem{example}[theorem]{Example}

\theoremstyle{remark}
\newtheorem{remark}[theorem]{Remark}

\newcommand{\ben}{\begin{enumerate}}
\newcommand{\een}{\end{enumerate}}
\newcommand{\bit}{\begin{itemize}}
\newcommand{\eit}{\end{itemize}}

\def\R {{\mathbb R}}

\def\N{{\mathbb N}}
\def\T{{\mathbb T}}

\def\Z {{\mathbb Z}}

\def\Homeo{{\mathrm{Homeo}}\,}

\def\d{{\delta}}
 
\def\eps{{\varepsilon}}

\def\F{{\mathcal F}}
\def\M{{\mathcal M}}

\def\Lip{\mathrm{Lip}} 
\def\Is{\mathrm{Is}} 
\def\Mol{\mathrm{Mol}} 

\def\a {\alpha}
\def\s {\sigma}
\def\g {\gamma}
\def\t {\tau}
\def\sk {\vskip 0.3cm}

\def\Is {\mathrm{Is}}
\def\0{\mathbf{0}}

\def\Unif{\operatorname{Unif}}

\def\RUC{\operatorname{RUC}}
\def\WAP{\operatorname{WAP}}

\begin{document}

\title[]{Lipschitz-Free Spaces: A Topometric Approach and Group Actions}    
 
		\author[]{Michael Megrelishvili}
		\address[Michael Megrelishvili]
	{\hfill\break Department of Mathematics
		\hfill\break
		Bar-Ilan University, 52900 Ramat-Gan	
		\hfill\break
		Israel}
		\email{megereli@math.biu.ac.il}

 \date{September 2025}    


\thanks{This research was supported by the Gelbart Research Institute at the Department of Mathematics, Bar-Ilan  University} 

\begin{abstract}  
 We introduce a topometric version of Lipschitz-free spaces and study its universal property. Another aim of this paper is 	
to investigate actions of topological groups $G$ on Lipschitz-free spaces $\mathcal{F}(M)$, induced by isometric actions on pointed metric spaces $M$. In particular, we study the associated dynamical $G$-systems under the weak-star topology, focusing on the dual action on $\mathrm{Lip}_0(M) = \mathcal{F}(M)^*$ and the bidual $\mathcal{F}(M)^{**}$.  
\end{abstract}

\maketitle
	
\setcounter{tocdepth}{1}
	\tableofcontents 

\section{Introduction}

The \textit{Lipschitz-free space} $\mathcal{F}(M)$ is a Banach space canonically associated with every pointed metric space $M$, serving as a powerful tool for the study of various metric properties. The theory of Lipschitz-free spaces has developed into a rapidly growing and significant area of research; see, for instance, \cite{Weaver, GK, Gode15, CDW, GPPR, Petit-Diss, APPP20, CDT, APS24} and the references therein.

This construction is also known by other names in the literature.  It is referred to as the \emph{Arens–Eells embedding}, following the foundational work in \cite{AE}, and as the \emph{free Banach space} over $M$, as in Pestov’s work \cite{Pe-free}. A variant of this construction already appears in classical optimization theory, particularly in the context of transportation problems. As a result, the associated norm is also known as the \emph{transportation cost norm}  \cite{Ost2-Med, Ost2-Anal}, the \emph{Kantorovich–Rubinstein norm} \cite{MPV}, or the \emph{Kantorovich norm} \cite{Ver}.       

In Section~\ref{s:defin}, we recall classical definitions and properties of Lipschitz-free spaces over pointed metric spaces $(M,d,\mathbf{0})$. 
In Section~\ref{s:Topom}, we introduce a \emph{topometric version} of Lipschitz-free spaces for pointed topometric spaces $\mathcal{M} := (M, d, \tau, \mathbf{0})$. Topometric spaces, introduced by Ben Yaacov \cite{BenTopom}, have proved important in recent developments; see, for example, \cite{BenTopom08, BenTopom, BBM, BMT17, Zuck}. Lower semicontinuous metrics compatible with topology arise naturally in many contexts, especially in Banach space theory, where the interplay between the norm and the weak-star (or weak) topology is central. 
We therefore expect that topometric generalizations of Lipschitz-free spaces offer a promising direction for future research. For the topometric version we have a continuous pairing $\F(\M) \times \Lip_0(M,d,\t) \to \R$, where $\Lip_0(M,d,\t)$ is a normed subspace of the dual $\F(\M)^*$. 
More precisely, by Theorem \ref{t:topomDual}, 
$\Lip_0(M,d,\t)=$ $\{f \in \mathcal{F}(\mathcal{M})^*: \ f|_M \in C(M,\t)\}$.  

On the other hand, $\Lip_0(M,d,\t)$ sometimes is quite computable and might be even separable (in contrast to the classical setting, where typically $\F(\M)$ contains an isomorphic copy of $\ell_1$). Note that $\Lip_0(M,d,\t)$ is complete for compact $\t$ and bounded $d$ (Proposition \ref{p:compBound}) and the boundedness is essential (see Example \ref{e:ForTopometricRevised}). 

For every isometric action of $G$ on a pointed metric space $(M,d,\0)$ which fixes $\0$ we have   
induced linear actions on the Lipschitz-free space $\mathcal{F}(M)$, on its dual, and bidual, primarily using the weak-star topology. 

Section~\ref{s:actions} is devoted to continuity aspects of natural induced dual actions under the weak-star topology. A recent result \cite[Proposition 2.3]{APS24} implies that the so-called \emph{Lipschitz realcompactification} $M^{\mathcal{R}}$ \cite{GarMer} of $(M,d)$ can be naturally identified with the weak-star closure $\overline{\delta(M)}^{w^*} \subset (\mathcal{F}(M))^{**}$ of $M$ in the bidual. 
In Theorem~\ref{t:positive} and Corollary~\ref{c:TildeCont}, we investigate when the canonically defined 
action of $G$ on $M^{\mathcal{R}}$ is continuous.  

To every Banach space $(V, \|\cdot\|)$, one may associate several important objects, such as the topological group $\mathrm{Is}_{\mathrm{lin}}(V)$ of all linear onto isometries (with the strong operator topology) and its canonical dual action on the weak-star compact unit ball $B_{V^*}$ of the dual space $V^*$. A natural approach is to linearize abstract continuous actions $G \times X \to X$ of a topological group $G$ on a topological space $X$ (i.e., a $G$-\textit{space}) via the dual action on $B_{V^*}$ for some $V$.    

\begin{definition} \label{d:repr} \cite{Me-FRepres03,GM-survey14}  
	Let $X$ be a $G$-space. A \emph{representation} of $(G,X)$ on a Banach space $(V,\|\cdot\|)$ is a pair
	\[
	h\colon G \to \mathrm{Is}_{\mathrm{lin}}(V), \quad \alpha\colon X \to V^*,
	\]
	where $h$ is a continuous homomorphism and $\alpha$ is a weak-star continuous bounded $G$-mapping (e.g., $\alpha(X) \subset B_{V^*}$) with respect to the dual action 
	\[
	G \times V^* \to V^*, \quad (g \varphi)(v) := \varphi(g^{-1}v)
	\]
	which makes the following diagram commutative (that is, $\a(gx)=g \a(x)$ for every $(g,x) \in G \times X$). 
	\begin{equation*} \label{diag2}
		\xymatrix{
			G \ar@<-2ex>[d]_{h} \times X \ar@<2ex>[d]^{\alpha} \ar[r] & X \ar[d]^{\alpha} \\
			\mathrm{Is}_{\mathrm{lin}}(V) \times V^* \ar[r] & V^*
		}
	\end{equation*}   
	A representation is called \emph{proper} if $\alpha$ is a topological embedding. Note that if $X$ is compact, then every weak-star continuous map $\alpha\colon X \to V^*$ is necessarily bounded.
\end{definition}

This framework provides new tools for analyzing abstract dynamical $G$-systems using the geometry of Banach spaces. For further background and applications, see \cite{GM-AffComp13, GM-survey14, Gl-YomDin, Me-b}.

In Section \ref{s:Repr} we propose studying representations of group actions on the Lipschitz-free space $V := \mathcal{F}(M)$, where $M$ is a pointed metric space and $\mathcal{F}(M)^* = \Lip_0(M)$. We focus on the case where $h$ is a homomorphism from $G$ into $\mathrm{Is}(M) \subset \mathrm{Is}_{\mathrm{lin}}(\mathcal{F}(M))$; see Definition~\ref{d:RepresFL}. 
Theorem~\ref{t:suff} confirms that there exist sufficiently many representations of compact $G$-spaces on Lipschitz-free spaces $\mathcal{F}(M)$.

Section \ref{s:horo}, and in particular Theorem~\ref{t:hori},  provides a particular instance of a representation on free Lipschitz spaces for isometric actions, involving the so-called \textit{metric (horo)compactification}
\[
\mu \colon M \to \widetilde{M} \subset \Lip_0(M),
\]
where $M$ is a pointed metric $G$-space with an isometric action. For bounded metric $G$-spaces, we derive further consequences for the \textit{Gromov $G$-compactification} (Definition~\ref{d:GromCom}). 

A standard definition of amenability for topological groups $G$ is the existence of a fixed point in every affine compact $G$-flow; see \cite[Theorem III.3.1]{Gl-book} and \cite[Theorem G.1.7]{BHV}. In the context of Lipschitz-free spaces and metric geometry, it is natural to consider weaker forms of amenability with a metric flavor. See Question~\ref{q:ind}.b and the 
discussion in Section~\ref{s:horo}. 

Theorem~\ref{t:WAP} states that for every isometric $G$-space $(M, d, \mathbf{0})$ with a bounded \textit{stable} metric $d$, the corresponding metric $G$-compactification $\widetilde{M}$ is a weakly almost periodic $G$-flow. By a Ryll–Nardzewski-type fixed point theorem (see Fact~\ref{f:measure}), it follows that for every $a \in M$, the internal metric functional
\[
\mu_a \in \Lip_0(M), \quad \mu_a(x) := d(a, x) - d(a, \mathbf{0}),
\]
is \emph{amenable} in the sense that the corresponding “cyclic” affine $G$-flow $\overline{\operatorname{co}}^{w^*}(G \mu_a)$ admits a $G$-fixed point. In fact, this property holds for every Lipschitz map $f \in \overline{\operatorname{co}}^{w^*}(\widetilde{M})$, where $\widetilde{M} \subset \Lip_0(M)$ denotes the metric (horofunction) compactification of the stable metric space $M$, represented in the Lipschitz-free space.

This applies, for example, when $M := B_V$ is the unit ball of a Banach space $(V, \|\cdot\|)$ with a \emph{stable} norm in the sense of Krivine–Maurey \cite{KrMa}, and $G$ is any subgroup of the linear isometry group $\mathrm{Is}_{\mathrm{lin}}(V)$. See Corollary~\ref{c:StableBanach} for further details.

Throughout this paper, we pose several questions: \ref{q:topometric}, \ref{q:SEPARABtopom}, \ref{q:light}, \ref{q:delta-starComp}, \ref{q:REPR}, \ref{q:ind}, \ref{q:WAP}, \ref{q:tame}, \ref{q:MetrAmen}.

\smallskip 
\noindent\textbf{Acknowledgment.}  
I am grateful to Marek Cuth, Michal Doucha and Colin Petitjean for several valuable suggestions. I also thank the organizers and participants of the \emph{First Conference on Lipschitz-Free Spaces} (Besançon, September 2023) for their inspiration and stimulating discussions. Part of the present article was presented at that conference. 

Special thanks we address to the referee for his/her very useful numerous concrete suggestions and constructive criticism which helped us essentially improve the present work. 

\sk 
\section{Lipschitz-free Banach spaces} 
\label{s:defin}

In this section, we briefly recall some well-known facts about Lipschitz-free spaces. Let $M$ be a nonempty set. A \textit{molecule} on $M$ is a formal finite sum 
\[ m = \sum_{i=1}^n c_i(x_i - y_i), \]
where $x_i, y_i \in M$, $c_i \in \mathbb{R}$, and $n \in \mathbb{N}$. Such a molecule can also be viewed as a finitely supported function $m \colon M \to \mathbb{R}$ satisfying $\sum_{x \in M} m(x) = 0$. The set $\Mol(M)$ of all molecules is a real vector space. 
Let $d$ be a pseudometric on $M$. Define a seminorm on $\Mol(M)$ by
\[ \|m\|_d := \inf \Bigl\{ \sum_{i=1}^n |c_i| d(x_i, y_i) : m = \sum_{i=1}^n c_i(x_i - y_i) \Bigr\} \]
over all finite decompositions. 
It is well known that $\|\cdot\|_d$ is a norm if and only if $d$ is a metric. In that case, $(\Mol(M), \|\cdot\|_d)$ is referred to as the \textit{Arens--Eells normed space} of $(M,d)$. We will usually write simply $\Mol(M)$ and $\|\cdot\|$. 
This norm is sometimes called the \textit{Kantorovich--Rubinstein norm} and plays a central role in optimization theory \cite{Ver, MPV, Ost2-Med, Ost2-Anal}.

The \textit{Lipschitz-free space} $\mathcal{F}(M)$ over $(M,d)$ is the completion of $(\Mol(M), \|\cdot\|)$, forming a Banach space.

Let $(M, d, \mathbf{0})$ be a pointed metric space. For each $x \in M$, define the molecule $\delta_x := x - \mathbf{0}$. The set $\{\delta_x : x \in M \setminus \{\mathbf{0}\}\}$ forms a Hamel basis for $\Mol(M, d)$. Define the canonical injection
\[ \delta \colon M \to \Mol(M), \quad x \mapsto \delta_x. \]
Clearly, $\delta_{\mathbf{0}}$ is the zero element of $\mathcal{F}(M)$. 
Any molecule can be expressed as $m = \sum_{i=1}^n r_i x_i = \sum_{i=1}^n r_i (\delta_{x_i} - \mathbf{0})$.  

Now consider the dual space $\mathcal{F}(M)^*$. For each functional $F \colon \mathcal{F}(M) \to \mathbb{R}$, we define a corresponding function $f \colon M \to \mathbb{R}$ by $f(x) := F(\delta_x)$. Conversely, given a function $f \colon M \to \mathbb{R}$ with $f(\mathbf{0}) = 0$, we extend it linearly to $F \colon \Mol(M) \to \mathbb{R}$ via
\[ F(m) = \sum_{i=1}^n c_i(f(x_i) - f(y_i)) \quad \text{for } m = \sum_{i=1}^n c_i(x_i - y_i). \]

Let $\Lip_0(M)$ denote the Banach space of all Lipschitz functions $f \colon M \to \mathbb{R}$ with $f(\mathbf{0}) = 0$, equipped with the norm $\|f\|_{\Lip} := \Lip(f)$.

\begin{fact} \label{f:facts} 
	Let $(M, d, \mathbf{0})$ be a pointed metric space. Then:
	\begin{enumerate}
		\item $\mathcal{F}(M)^* = \Lip_0(M)$.
		\item The weak-star topology on bounded subsets of $\Lip_0(M)$ coincides with the topology of pointwise convergence.
		\item The norm $\|\cdot\|_d$ on $\Mol(M)$ is the largest seminorm satisfying $\|\delta_x - \delta_y\| \le d(x, y)$, and in fact $\|\delta_x - \delta_y\| = d(x, y)$. Thus, $\delta \colon M \to \Mol(M)$ is an isometric embedding.
		\item (Universal property) Let $V$ be a Banach space and $f \in \Lip_0(M, V)$. There exists a unique linear map $T_f \in L(\mathcal{F}(M), V)$ such that $f = T_f \circ \delta$ and $\|T_f\| = \|f\|_{\Lip}$.
		\item (Canonical linearization) For any Lipschitz map $f \colon (M_1, \mathbf{0}) \to (M_2, \mathbf{0})$ between pointed metric spaces, there exists a unique continuous linear map $\bar{f} \colon \mathcal{F}(M_1) \to \mathcal{F}(M_2)$ such that $\bar{f} \circ \delta_1 = \delta_2 \circ f$ and $\|\bar{f}\| = \|f\|_{\Lip}$.
	\end{enumerate}
\end{fact}

\begin{remark} \label{r:NORM} 
	Let $V$ be a Banach space and $i \colon V \hookrightarrow V^{**}$ its canonical isometric embedding into the bidual. In particular, this holds for $\mathcal{F}(M)$. 
	Thus, we have the chain of isometric embeddings:
	\[ M \xrightarrow{\delta} \mathcal{F}(M) \xrightarrow{i} \mathcal{F}(M)^{**} = \Lip_0(M)^*. \]
	We continue to denote the composition $i \circ \delta$ simply by $\delta$, identifying $M$ with its image $\delta(M)$. Then, for every $v \in \mathcal{F}(M)$,
	\[ \|v\| = \sup \{ f(v) : f \in \Lip_0(M), \ \|f\|_{\Lip} \le 1 \}. \]
	The space $\mathcal{F}(M)$ can thus be defined as
	\[ \mathcal{F}(M) = \overline{\operatorname{span}\{\delta(M)\}}^{\|\cdot\|}. \]
\end{remark}

\begin{remark} \label{r:w=w-st} 
	Let $V$ be a Banach space. Then the weak topology on $V$ coincides with the weak-star topology on $V \subset V^{**}$ inherited via the canonical embedding $i \colon V \hookrightarrow V^{**} = (V^*)^*$.
\end{remark}

\begin{proposition} \label{p:w=n} 
	Let $(M, d, \mathbf{0})$ be a pointed metric space.
	\begin{enumerate}
		\item The map $\delta \colon M \to (\mathcal{F}(M), \text{weak})$ is a topological embedding. Hence, the weak and norm topologies coincide on $\delta(M) \subset \mathcal{F}(M)$.
		\item The map $\delta \colon M \to (\Lip_0(M)^*, \text{weak}^*)$ is a topological embedding. Hence, the weak-star and norm topologies coincide on $\delta(M) \subset \mathcal{F}(M)^{**}$.
	\end{enumerate} 
\end{proposition}

\begin{proof}
	(1) Suppose a net $(m_i)$ in $M$ converges weakly to $a \in M$. Then for every $f \in \mathcal{F}(M)^*$, we have $\lim f(m_i) = f(a)$. In particular, take $f = \mu_a$ with $\mu_a(x) := d(a, x) - d(a, \mathbf{0})$. One can verify that $\mu_a \in \Lip_0(M)$ (see Theorem~\ref{t:hori}). On the other hand,
	\[ \mu_a(m_i) - \mu_a(a) = d(m_i, a) = \|m_i - a\|. \]
	Hence, weak convergence implies norm convergence, and the topologies coincide.
	
	(2) Apply Remark~\ref{r:w=w-st} to $V = \mathcal{F}(M)$ and use part (1) together with Remark~\ref{r:NORM}, identifying $M = \delta(M) \subset \mathcal{F}(M)^{**} = \Lip_0(M)^*$.
\end{proof}

\smallskip

\noindent
Proposition~\ref{p:w=n}(1) is well known; see, for example, \cite[Lemma 1.2.3]{Petit-Diss}. The completeness assumption on $M$ is not essential at this point. In fact, if $M$ is complete, then $M$ is weakly closed in $\mathcal{F}(M)$.

\sk  
\section{Topometric spaces and their Lipschitz-free spaces}
\label{s:Topom}

According to Proposition~\ref{p:w=n}, the topological space $M$ can be identified with $\delta(M)$ equipped with the weak-star topology inherited from $\Lip_0(M)^*$. The metric induced by the norm on the weak-star closure $cl_{w^*}(M) = \overline{M}^{w^*}$ is lower semicontinuous with respect to the weak-star topology. 
Moreover, if $(M, d)$ is bounded, then its isometric image $\delta(M)$ is norm-bounded, and $\overline{M}^{w^*} \subset \Lip_0(M)^*$ is weak-star compact.

On the other hand, by a well-known result of Jayne, Namioka, and Rogers \cite[Theorem~2.1]{JNR90}, every compact space $(K, \tau)$ with a bounded lower semicontinuous metric $d$ can be represented in the dual of a Banach space $V^*$ such that the norm induces the metric $d$ on $K$, and the weak-star topology on $V^*$ induces the original topology $\tau$. A significantly simpler proof was later given by Raja \cite[Theorem~2.3]{Raja24}, who also noted a resemblance to the theory of Lipschitz-free spaces.

In Theorem~\ref{t:topomDual}, we extend this result to a more general setting that of not necessarily compact spaces under natural assumptions. Specifically, for \emph{completely regular topometric spaces}. 

First we recall necessary definitions. 
Recall that a metric $d \colon X^2 \to [0,\infty)$ is said to be lower semicontinuous on a topological space $(X,\tau)$ if the set $\{(x,y): X^2: d(x,y) \leq c\}$ is closed in $(X,\tau)\times (X,\tau)$ for every positive $c$. 

A \emph{topometric space}, as introduced by I. Ben Yaacov \cite{BenTopom08, BenTopom}, is a triple 
$\mathcal{M} := (M, d, \tau)$, where $d$ is a lower semicontinuous metric on $(M,\tau)$ that refines the topology $\tau$ on $M$. 
In fact, the original definition in \cite{BenTopom} was for the $[0,\infty]$-valued metrics.  
In contrast, we do not allow in the present paper metrics with infinite values. Thus, $d(x,y)$ always will mean a nonnegative real number. 

 The notion of topometric spaces has found numerous important applications; see, for example, \cite{BenTopom08, BenTopom, BBM, BMT17, Zuck}.
One of the original motivations for the concept  arose from the space of types in first-order logic.  
Moreover, 
any subset $M$ of a dual Banach space $V^*$, equipped with the norm-induced metric and the weak-star topology, provides a natural example of a topometric space. Conversely, Theorem~\ref{t:topomDual} below shows that many important topometric spaces can, in fact, be realized as subsets of dual Banach spaces.

Assume that we are given a \emph{pointed} topometric space  
$\mathcal{M} := (M, d, \tau, \0)$. Our goal is to examine a topometric generalization of Lipschitz-free spaces.   
We follow an approach similar to the traditional definition  described in Remark~\ref{r:NORM}. For the topometric space $\mathcal{M}$, consider 
\[
V := \Lip_0(M, d) \cap C(M, \tau),
\]
as a normed subspace of $(\Lip_0(M), \|\cdot\|_{\Lip})$. For simplicity, we denote it by $V$ or $\Lip_0(M, d, \tau)$.  

In the dual Banach space $V^*$, for every $x \in M$, define the evaluation functional $\delta_x \colon V \to \mathbb{R}$ by $\delta_x(f) := f(x)$ for every $f \in V$. Clearly, $\delta_x$ is linear and $\d(\0)=0$. It is also $\|\cdot\|_{\Lip}$-continuous on $V$, as shown by the following inequality, valid for all $f_1, f_2 \in V$:
\[
|\delta_x(f_1) - \delta_x(f_2)| = |f_1(x) - f_2(x)| = |(f_1 - f_2)(x) - (f_1 - f_2)(\0)| \leq \ d(x, \0) \cdot \|f_1 - f_2\|_{\Lip}.
\]
Thus, $\|\delta_x\| \leq d(x, \0)$. Since $d(x,0)$ is fixed for any given $x \in M$, $\d_x$ is $||\cdot||_{\Lip}$-continuous. Thus, indeed $\d_x \in V^*$ and the following function is well defined:  
$$
\d \colon M \longrightarrow V^*=\Lip_0 (M,d,\t)^*, \ \ \d(x):=\d_x, \ \ \d_x(f)=f(x) \quad(f\in V). 
$$ 

 \begin{definition} \label{d:topomLF}  
 	Let $\mathcal{M} = (M, d, \tau, \mathbf{0})$ be a pointed topometric space. The \emph{\textbf{topometric Lipschitz-free space}} associated with $\mathcal{M}$ is defined as the norm-closed linear span
 	\begin{equation} \label{eq:topo-alt} 
 		\mathcal{F}(\mathcal{M}) := 
 		\overline{\mathrm{span} \{\delta(M)\}}^{\|\cdot\|} \subseteq \Lip_0(M,d,\tau)^*,
 	\end{equation}
 	where $\Lip_0(M,d,\tau)$ denotes the space $V := \Lip_0(M,d) \cap C(M,\tau)$ equipped with the Lipschitz norm $\|\cdot\|_{\Lip}$. 
 	Thus, for every $v \in \mathcal{F}(\mathcal{M})$, the norm is given by
 	\begin{equation} \label{eq:NewNorm} 
 	\|v\| = \sup\left\{ \langle f, v \rangle : f \in \Lip_0(M,d,\tau), \ \|f\|_{\Lip} \leq 1 \right\}. 	
 	\end{equation}
 	An alternative name for $\F(\M)$ might be \emph{the topometric Arens-Eells space}. 
 	
 	In the special case where the topology $\tau$ coincides with the metric topology induced by $d$, we recover the classical Lipschitz-free space $\mathcal{F}(M,d)$, as described in Remark~\ref{r:NORM}. 
 \end{definition}
 
 
 \smallskip
 
 \noindent
 Recall the notion of \textit{completely regular} topometric spaces, introduced by Ben Yaacov~\cite{BenTopom}.
 
 \begin{definition} [\cite{BenTopom}] \label{d:normal} 
 	A topometric space $\mathcal{M} = (M, d, \tau)$ is said to be \textbf{completely regular} if the family of all $\tau$-continuous 1-Lipschitz functions
 	\[
 	C_{L1}(M) := \left\{ f \in C(M,\tau) : |f(x) - f(y)| \leq d(x,y) \ \text{for all } x, y \in M \right\}
 	\]
 	is \emph{sufficient} for $\mathcal{M}$, in the sense that the following two conditions are satisfied:
 	\begin{enumerate}
 		\item For every point $x_0 \in M$ and every $\tau$-closed subset $F \subset M$ with $x_0 \notin F$, there exists a function $f \in C_{L1}(M)$ and distinct real numbers $a \ne b$ such that $f(x_0) = a$ and $f(F) = b$.
 		
 		\item For every pair $x_1, x_2 \in M$, the metric $d$ satisfies
 		\[
 		d(x_1, x_2) = \sup\left\{ |f(x_1) - f(x_2)| : f \in C_{L1}(M) \right\}.
 		\]
 	\end{enumerate}
 \end{definition}

\begin{lemma} \label{l:still}  
	In terms of Definition \ref{d:normal} define a ``pointed version" of $C_{L1}(M)$ as
	\[
	C_{L1}(M,\mathbf{0}) := \left\{ f-f(\0): \  f \in C_{L1}(M)\right\}.  
	\] 
	Then 
	\begin{equation} \label{eq:substract} 
	C_{L1}(M,\mathbf{0})  = \left\{f \in C_{L1}(M): \  |f(x) - f(y)| \leq d(x,y) \ \ \text{for all} \ x,y \in M, \ \text{and} \ f(\mathbf{0}) = 0 \right\} 	
	\end{equation} 
	 holds and the family $C_{L1}(M,\mathbf{0})$ is still \textbf{sufficient}.
\end{lemma}  
\begin{proof} 
	Indeed, take into account that the Lipschitz constants of $f$ and $f-f(\0)$ are the same and  $|(f(x_1)-f(\0))-(f(x_2)-f(\0))|=|f(x_1)-f(x_2)|$. 
	Also, 
	if $f(x_0)=a \neq b = f(F)$ for some $f \in C_{L1}(M)$, then  $f-f(\0) \in C_{L1}(M,\0))$ separates $x_0$ and $F$.    
\end{proof}

\begin{remark} \label{r:TopomExamples}
	Completely regular topometric spaces form a rich and useful class that is closed under taking subspaces. We list below several notable examples, as discussed in~\cite{BenTopom,BBM,BMT17,Zuck}:
	\begin{enumerate}
		\item (\emph{Classical case}) $(M,d)$ is a metric space and $\tau$ is the topology induced by $d$.
		
		\item Every subset of a normed space  equipped with the norm metric and the weak topology.
		
		\item Every subset in the dual of a normed space with its dual norm and the weak-star topology.
		
		\item The triplet $(S(M), \partial, \tau)$, where $M \to (S(M),\t)$ is the Samuel compactification of a bounded metric space $(M,d)$ and
		\[
		\partial(u,v) := \sup\{ |f(u) - f(v)| : f \in C_{L1}(M) \}.
		\]
	
		\item Topometric spaces $(K,d,\tau)$ where $(K,\t)$ is compact (or more generally, any \emph{normal topometric space} in the sense of~\cite{BenTopom}).
		
		\item (\emph{Topometric groups}) The triplet $(G,d_u,\t)$, where $(G,\tau)$ is a metrizable topological group, $d_L$ is a compatible left-invariant bounded metric on $(G,\tau)$ and  
		\[
		d_u(g_1,g_2) := \sup \{ d_L(g_1 h, g_2 h) : h \in G \}.
		\]  
	\end{enumerate}
\end{remark}

\sk  
\begin{theorem} \label{t:topomDual} 
Let $\mathcal{M}:=(M,d,\t,\0)$  be a pointed completely regular 
(in the sense of Definition \ref{d:normal}) 
topometric space. 
Then 
\begin{enumerate}
	\item The inherited metric and the weak-star topology on the subset $\d(M)$ of the dual $\Lip_0(M,d,\t)^{*}$ gives the original topometric structure on $(M,d,\t)$. 
	\begin{enumerate}
		\item $\d \colon (M,d) \to \mathcal{F}(\mathcal{M}) \subset (\Lip_0(M,d,\t)^{*},||\cdot||)$ is an isometric embedding. 
			\item 
		$\d \colon (M,\t) \to \mathcal{F}(\mathcal{M}) \subset (\Lip_0(M,d,\t)^{*},weak^*)$ is a topological embedding. 
	\end{enumerate}
	
	\item $\{\d(x): x \in M \setminus \{\0\}\}$ is linearly independent in $\mathcal{F}(\mathcal{M})$.   
		\item Define \ 
		$
		\mathcal{F}(\mathcal{M})_{\t}^{*}:=\{f \in \mathcal{F}(\mathcal{M})^*: f|_M \in C(M,\t)\}.	
		$ Then we have a natural isometric isomorphism 
		\begin{equation} \label{eq:form} 
			T \colon \mathcal{F}(\mathcal{M})_{\t}^{*} \to \Lip_0 (M,d,\t), \ \ \ T(\psi)=\psi \circ \d. 	
		\end{equation} 	
\end{enumerate} 
\end{theorem}
\begin{proof} (1)  
We proceed mainly as in the proofs of \cite[Theorem 2.3]{Raja24} and \cite[Theorem 2.1]{JNR90}, where the topology $\t$ was compact and $d$ is bounded.  

Observe that $C_{L1}(M,\0)$ (from Lemma \ref{l:still}) is exactly the closed unit ball $B_V$  of $V:=\Lip_0 (M,d,\t)$. Thus, by Definition \ref{d:topomLF} (for the vector $v:=\d(x)-\d(y)$),  and the sufficiency condition (2) for $C_{L1}(M,\0)$ (Definition \ref{d:normal} and Lemma \ref{l:still} ), we obtain that 
$$
d(x,y) = 
\sup \{|f(x)-f(y)|: \ f \in C_{L1}(M,\0)\} = ||\d(x)-\d(y)||_{V^*}. 
$$ 
This proves (a). 
\sk  
(b) 
Let $(x)_i$ be a converging in $(M,\tau)$ net with $\lim (x_i)=x$. 
In order to check the weak-star continuity of $\d \colon (M,\t) \to (\Lip_0(M,d,\t)^{*},weak^*)$, we have to establish that $\lim (\delta(x_i))=\delta(x)$. 
Since $\langle \delta(a),f \rangle=f(a)$ for every $a \in M$ and $f \in \Lip_0(M,d,\t)$,  
it is equivalent to show that $\lim (f(x_i))=f(x)$.
However, this is clear because $f \in C(M,\t)$ is $\t$-continuous.

Recall that for every normed space $V$ the dual space $(V^*,weak^*)$ in its weak-star topology naturally is embedded topologically into the power (product) $\R^{B_V}$. 
More precisely, the diagonal map 
$$\Delta \colon \delta(M) \to \R^{B_V}, \ \Delta(\delta(x)):=(f(x))_{f \in B_V}$$   
is a natural weak-star  embedding. Since, $\d \colon M \to \d (M)$ is weak-star continuous, the family of functions $\{f \circ \d \colon M \to \R: f \in B_V\}$ is continuous. The following diagonal product is a topological embedding because the family $B_V$ of functions separates points and closed subsets of $M$ by virtue of  the sufficiency condition (1) for $C_{L1}(M,\0)=B_V$ (Lemma \ref{l:still} ). 
This implies that 
$$\d \colon (M,\t) \to (\Lip_0 (M,d,\t)^*, weak^*)$$
is a topological embedding.  

\sk  
(2)  
Let $A:={x_1, \cdots, x_n} \subseteq M \setminus \{\0\}$ be a finite subset. For a given $1 \leq i \leq n$ define the finite set 
$F_i:=\{\0\} \cup (A \setminus \{x_i\})$.  Definition \ref{d:normal} and Lemma \ref{l:still} show that  there exist: distinct $a \neq b$ and a $1$-Lipschitz $\t$-continuous function $f_i \colon M \to \R$ from $B_V=C_{L1}(M,\0)$ such that  $f_i(F_i)=a \neq b = f_i(x_i)$. 
Since $\0 \in F_i$ and $f_i \in C_{L1}(M,\0)$, one may suppose, in addition, that $f_i(F_i)=a=0$. Thus, $f_i(x_i) \neq 0$.  
 This guarantees that $\d(x_i) \notin \operatorname{span} (\d(A \setminus \{x_i)\})$. So, $\d(A)$ is linearly independent in $\mathcal{F}(\mathcal{M})$. 
 
 \sk
 (3) 
 We verify that the following restriction assignment 
 (\ref{eq:form}) is an isometric isomorphism  
 $$T \colon \mathcal{F}(\mathcal{M})_{\t}^{*} \to \Lip_0 (M,d,\t), \ \ \ T(\psi)=\psi \circ \d. 
 $$ 	
 Clearly, 
 $$
 |T(\psi)(x)-T(\psi)(y)| \leq ||\psi||\cdot ||\d(x)-\d(y)||= ||\psi||\cdot d(x,y). 
 $$
 Thus, $||T(\psi)||_{\Lip} \leq ||\psi||$. Also, $T(\psi)(\0)=0$. 
 Therefore, $T(\psi) \in \Lip_0 (M,d)$. 
 
 
 Since, the bounded functional $\psi \colon \mathcal{F}(\mathcal{M}) \to \R$
 belongs to $\mathcal{F}(\mathcal{M})_{\t}^{*}$, the composition 
 $T(\psi)=\psi\circ \d$ is  $\t$-continuous. Thus, $T(\psi) \in C(M,\t)$.   
 This means that indeed $T(\psi) \in \Lip_0 (M,d,\t)$ and $T$ is well defined. The remainder of the arguments now follows exactly as in the classical setting. We give the details for the sake of completeness. 
 
 \sk 
 \noindent {\bf Claim 1:} $T \colon \mathcal{F}(\mathcal{M})_{\t}^{*} \to \Lip_0 (M,d,\t)$ is a linear isometric embedding.  
 \sk 
 For every $\psi \in \mathcal{F}(\mathcal{M})_{\t}^{*}$ and every $x \in M$ we have 
 $$
 \langle T(\psi),\d(x)\rangle=T(\psi)(x)=(\psi\circ \d)(x)=
 \langle \d(x),\psi\rangle.  
 $$
 Since $span (\d(M))$ is norm dense in $\mathcal{F}(\mathcal{M})_{\t}^{*}$, we obtain that 
 $$
 \langle T(\psi),v\rangle=\langle v,\psi \rangle \ \ \ \text{for all} \ v \in \mathcal{F}(\mathcal{M}).
 $$
 Using this equality we get 
 $$
 ||\psi||=\sup \{\langle v,\psi \rangle: v \in B_{\mathcal{F}(\mathcal{M})}\} = \sup \{\langle T(\psi),v \rangle: v \in B_{\mathcal{F}(\mathcal{M})}\} \leq 
 $$ 
 $$
 \leq \sup \{\langle T(\psi),w \rangle: w \in B_{\Lip_0 (M,d,\t)^*}\}=||T(\psi)||_{\Lip}.
 $$
 Thus, $||\psi| \leq ||T(\psi)||_{\Lip}.$ We already saw the converse direction $||T(\psi)||_{\Lip} \leq ||\psi||$. Therefore, we conclude that $||T(\psi)||_{\Lip} = ||\psi||$ for every $\psi \in \mathcal{F}(\mathcal{M})_{\t}^{*}$. 
 
 \sk 
 \noindent {\bf Claim 2:} $T \colon \mathcal{F}(\mathcal{M})_{\t}^{*} \to \Lip_0 (M,d,\t)$ is onto (hence, a linear isometry).   
 
 \sk 
 Let $f \in \Lip_0 (M,d,\t)$. Our aim is to find $\psi_f \in \mathcal{F}(\mathcal{M})_{\t}^{*}$ such that $T(\psi_f)=f$. Without any loss of generality we can assume that $||f||_{\Lip}=1$. 
 First we define the function $\psi_f$ on the subset $\d(M)$ by the rule
 $\psi_f(\d(x))=f(x)$ for every $x \in M$. Clearly, $\psi_f(0)=f(\0)=0$.

 By (2), 
 $\{\d(x): x \in M \setminus \{\0\}\}$ is linearly independent in $\mathcal{F}(\mathcal{M})$. Therefore, we canonically can extend this map to the linear map $\mathrm{span} \{\d(M)\} \to \R$ which we also denote by $\psi_f$. By Definition \ref{d:topomLF}, this map is continuous because for every $v \in \mathrm{span} \{\d(M)\}$ we have 
 $$|\psi_f(v)|=|f(v)|=|\langle f,v \rangle| \leq ||f||_{\Lip}\cdot ||v||=||v||.$$  	
 Since $span \{\d(M)\}$ is dense in $\mathcal{F}(\mathcal{M})$, there exists a unique linear continuous extension $\psi_f \colon \mathcal{F}(\mathcal{M}) \to \R$. 
 Furthermore, by construction $\psi_f \circ \d=f \in C(M,\t)$ holds. Therefore, 
 $\psi_f \in \mathcal{F}(\mathcal{M})_{\t}^{*}$ is such that $T(\psi_f)=f$.   
\end{proof}

Separability of the Banach space $\F(\M)$ is equivalent to separability of $(M,d)$ as it follows by Theorem \ref{t:topomDual}.1.a and Definition \ref{d:topomLF}. 

\begin{remark} \label{r:converseTOPOM} 
One of the consequences of Theorem \ref{t:topomDual} is that a topometric pointed space 
is completely regular if and only if it can be embedded (topometrically) into the dual Banach space $V^*$, where $V$ is a normed space. That is, example (3) from Remark \ref{r:TopomExamples} is universal, in a sense. 	
\end{remark}

\begin{example}[Compact topometric space with unbounded metric] \label{ex:NotBounded} 
	Let $K = \{0\} \cup \{\frac{1}{n}: n \in \N\}$ with its compact topology. Consider on $K$ the metric $d$ induced from the natural bijection $K \to \{0\} \cup \N$. That is, $d(\frac{1}{n},\frac{1}{m})=|n-m|$  and $d(0,\frac{1}{n})=n$ for every $n,m \in \N$.
	 Then $\mathcal{M}:=(K,d,\tau,0)$ is a pointed metric space with compact topology $\tau$ and an unbounded metric $d$.
	 \begin{proof}
 It is easy to see that $d$ is a lower semicontinuous metric with respect to $\tau$. Indeed, we show that 
for each \(c>0\) the set  
$U_c
=\{(x,y)\in K\times K: d(x,y)>c\}$ 
is open in $K \times K$.  Let $(x,y) \in U_c$, that is, $d(x,y) >c$. We check  that every \((x,y)\) admits an open neighbourhood $U \times V$ in \(K\times K\) contained in \(U_c\). 

\medskip
\noindent\textbf{Case 1.} $x=\frac{1}{n}$ and $y=\frac{1}{m}$ are both nonzero. Then take the singleton $U \times V:=\{(\frac{1}{n}, \frac{1}{m})\}$. 

\medskip
\noindent\textbf{Case 2.} Exactly one coordinate is \(0\); say
\((x,y)=(0,\frac{1}{n})\) with \(d(0,\frac{1}{n})=n > c\) (the other case is symmetric).
Set $V=\{\frac{1}{n}\}, \ U=\{0\}\,\cup\,\{\frac{1}{m}: m \geq 2n\}.$ 	 		
\end{proof}	 
\end{example} 

\sk  
Let \(V\) be a normed space (not necessarily complete) and \(\widetilde V\) its Banach completion.
One may identify \(V^*\) isometrically with \(\widetilde V{}^*\) by continuous extension of each functional from \(V\) to \(\widetilde V\). 
Then we have two comparable natural weakened (pointwise)  topologies on \(V^*\):
\begin{itemize}
	\item (weak\(^*\)‐topology) \(\sigma(V^*,V)\) is generated by the seminorms \(f\mapsto|f(v)|\), \(v\in V\).
	\item \(\sigma(V^*,\widetilde V)\) is generated by the seminorms \(f\mapsto|f(\hat v)|\), \(\hat v\in\widetilde V\).
\end{itemize}
Since \(V\subset\widetilde V\) (densely), clearly
\[
\sigma(V^*,V)\;\subseteq \;\sigma(V^*,\widetilde V),
\]
and in general the inclusion is strict when \(V\) is not complete.

\medskip
\noindent
However, on every norm‐bounded subset of \(V^*\) the two topologies coincide.  Equivalently,
\[
\text{a net }(f_\alpha)\subset V^*\text{ is weak\(^*\)‐convergent in }\sigma(V^*,\widetilde V)
\ \leftrightarrow \
(f_\alpha)\text{ is bounded and converges pointwise on} \ V.
\]

\begin{remark} \label{r:BanSteinh} 
Let $\mathcal{M}=(K,d,\tau,0)$	be the topometric space from Example \ref{ex:NotBounded}. Then by Theorem \ref{t:topomDual} we have a topometric embedding $\delta \colon \mathcal{M} \to V^{*}$, where $V:=\Lip_0(M,d,\t)$.  That is, $\d(K)$ is an unbounded metric subspace (isometric to $(K,d)$) of the dual Banach space $V^{*}$ and $\d(K)$  (being homeomorphic to $(K,\t)$) is \textbf{compact} with respect to the weak-star topology $w^*=\sigma(V^*,V)$ on $V^*$. Then, the pointwise topology $\sigma(V^*,\widetilde{V})$ restricted on $\d(K)$ is not compact and strictly stronger than the compact topology inherited from $\sigma(V^*,V)$. Indeed, otherwise we get a contradiction to one of the main corollaries to Banach-Steinhaus theorem: for every Banach space $W$, any weak-star compact subset $A$ in the dual $W^*$ is necessarily norm bounded (in the following setting: $W:=\widetilde{V}$ and $A:=\delta(K)$).   

This discussion has some consequences:
\begin{enumerate}
	\item  The inclusion $\Lip_0(M,d,\tau)=\mathcal{F}(\mathcal{M})_{\t}^{*} \subseteq \mathcal{F}(\mathcal{M})^*$ might be proper. Take for instance the compact topometric space $\mathcal{M}=(K,d,\tau,0)$ from Example \ref{ex:NotBounded} (compare Question \ref{q:topometric}). 
	\item The norm space $V:=\Lip_0(M,d,\t)$ (and, hence, also $\mathcal{F}(\mathcal{M})_{\t}^{*}$) might be not complete. For a concrete description see Example \ref{e:ForTopometricRevised}.   
	\item A gentle and quick reminder that Banach-Steinhaus theorem is not valid for normed spaces $V$. Namely, weak-star compact subsets $K$ in the dual $V^*$ might be norm unbounded.  
\end{enumerate}
\end{remark}

\begin{question} \label{q:topometric} 
	For which completely regular topometric spaces holds
	\[
	\mathcal{F}(\mathcal{M})^* = \Lip_0(M,d,\tau)\ ?
	\] 
\end{question}

By Definition \ref{d:topomLF}, $\mathcal{F}(\mathcal{M})$ is a Banach subspace of the dual space $\Lip_0 (M,d,\t)^*$. Consider the induced (continuous) bilinear map 
$$
w \colon \mathcal{F}(\mathcal{M}) \times \Lip_0 (M,d,\t) \to \R , \ \langle v,f \rangle:=f(v).
$$ 
Clearly, the normed space $\Lip_0(M,d,\tau)$ separates points of its dual. Hence, also points of $\mathcal{F}(\mathcal{M})$.  
Since the molecules separate the points of $\Lip_0 (M,d,\t)$, we obtain that $w$ separates the points on both sides. Therefore, we have a \textit{pairing}. 
It follows that $\Lip_0(M,d,\tau)$ is weak-star dense in $\mathcal{F}(\mathcal{M})^*$.

Since $\d \colon (M,d) \to \mathcal{F}(\mathcal{M})$ is isometric (by Theorem \ref{t:topomDual}.1.a), 
the corresponding norm on the molecules is 
``compatible" with $d$ in terms of \cite[Section 1.1]{MPV}.   
  
 The equality $\mathcal{F}(\mathcal{M})^* = \Lip_0(M,d,\tau)$ does not always hold (as we have seen in Remark \ref{r:BanSteinh}).
Even when $\Lip_0(M,d,\tau)$ is a Banach space (Example \ref{ex:10metric}), 
it need not exhaust the full dual of $\F(\M)$.  
Despite this, the equality 
 $\Lip_0(M,d,\t)=\mathcal{F}(\mathcal{M})_{\t}^{*}$, gives a chance for a fair computability. Moreover, 
 the presence of the natural continuous bilinear pairing $w$  provides a meaningful duality framework. Observe also, that sometimes $\Lip_0(M,d,\tau)$ are small, e.g. separable, as we  can see in concrete topometric samples of Examples  \ref{ex:10metric} and \ref{e:ForTopometricRevised}. While, in the classical case, 
 $\Lip_0(M,d)=\F(M)^*$ and, since 
 usually $\F(M)$ contains a copy of $\ell_1$, hence $\Lip_0(M,d)$ is not separable. This observation might be interesting for topometric $G$-spaces in the context of (metric) amenability; see the discussion after Fact \ref{f:measure}.

 \begin{question} \label{q:SEPARABtopom} 
 For which pointed topometric spaces $\M$ the normed space $\Lip_0(M,d,\tau)$ is separable~? 	Complete~?
 \end{question} 

 For some related results see Examples \ref{ex:10metric}, \ref{e:ForTopometricRevised}, Proposition \ref{p:compBound}, and Remark \ref{r:BALLS}. 

 
 \begin{remark} \label{r:Pet} 
 	Note that in \cite{GPPR} and \cite{APS24} several results 
 	deal with lower semicontinuous metrics and spaces of the form $\Lip_0(M, d) \cap C(M, \tau)$ (attributes of topometric  spaces). See, for example, 
 	\cite[Section 3 and Corollary 3.9]{GPPR} (where the authors study preduals of $\F(M)$) and \cite[Proposition 2.4 and Theorem 4.3]{APS24}.
 	These results probably are related to Question  \ref{q:topometric}. 
 	 I am grateful to M. Cuth and C. Petitjean for pointing this out to me.
 \end{remark}  
 
 \begin{corollary} \label{c:w=n} 
 	Let $\M:=(M,d,\t,\0)$ be a pointed completely regular topometric space. Then the corresponding pointwise topology 
 	$\sigma(\Lip_0(\mathcal{M}),\mathcal{F}(\mathcal{M}))$ on $\Lip_0(\mathcal{M})$ (induced by the pairing) coincides on norm \underline{bounded} sets with the topology of $M$-pointwise convergence.  
 \end{corollary}
 \begin{proof} 
 	By Theorem \ref{t:topomDual},  $\Lip_0 (M,d,\t) =\mathcal{F}(\mathcal{M})_{\t}^{*} \subseteq \mathcal{F}(\mathcal{M})^{*}$ holds.  
 	Let$(f_i)$ be a net in $\Lip_0 (M,d,\t)$. Then $\lim f_i = f$ in $\sigma(\Lip_0(\mathcal{M}),\mathcal{F}(\mathcal{M}))$-topology (or, the $w^*$-star topology inherited from the dual $\mathcal{F}(\mathcal{M})^{*}$)  means that $\langle u, f_i \rangle \to \langle u, f \rangle$ for every $u \in \mathcal{F}(\mathcal{M})$. 
 	Assume, in addition, that every $f_i$ (and, hence, also $f$)  belongs to a normed bounded subset. Then, since $\mathrm{span} (\d(M))$ is norm dense in $\mathcal{F}(\mathcal{M})$, $\lim f_i = f$ iff $\langle x, f_i \rangle \to \langle x, f \rangle$ for every $x \in M.$	
 \end{proof} 
 

 \begin{example} \label{e:measures}
 	Let \( K \) be a compact Hausdorff space and consider the Banach space \( C(K) \) with the supremum norm \( \|\cdot\|_{\sup} \). Then the closed unit ball \( \mathcal{M} := (B_{C(K)}, d, \tau_w) \), where \( d \) is the norm metric and \( \tau_w \) is the weak topology, is a completely regular topometric space (see Remark~\ref{r:TopomExamples}.2).
 	
 	In this case, there is a natural \textbf{injective} norm continuous linear map
 	\[
 	j \colon C(K)^* \to \Lip_0(B_{C(K)}, d, \tau_w), \quad j(\mu)(x) := \mu(x),
 	\]
 	for all \( \mu \in C(K)^* \) and \( x \in B_{C(K)} \). Each \( j(\mu) \) is affine and weakly continuous, and satisfies
 	\[
 	|j(\mu)(x_1) - j(\mu)(x_2)| = |\mu(x_1 - x_2)| \leq \|\mu\| \cdot \|x_1 - x_2\| = \|\mu\| \cdot d(x_1, x_2),
 	\]
 	so \( j(\mu) \) is \( \|\mu\| \)-Lipschitz. Hence, \( j \) is a linear operator of norm \( \leq 1 \), and its image belongs to  \( \Lip_0(B_{C(K)}, d, \tau_w) \).
 	
 	An additional benefit of this example is that for every molecule \( x - y \in \mathcal{M} = B_{C(K)} \), the Hahn–Banach theorem ensures that there exists a regular measure \( \mu \in C(K)^* \), with \( \|\mu\| = 1 \), such that \( j(\mu)(x - y) = d(x, y) \). In this setting, the supremum in the norm formula from Definition~\ref{d:topomLF} is always attained. 
 \end{example}

 Another case when the norm formula in  (\ref{eq:NewNorm}) is attaining is for \textit{compact topometric spaces} $(M,d,\t)$ (that is, with compact $\t$). This follows from a result of Matouskova \cite[Corollary 2.5]{Mat}. 
 
 It would be interesting to study the norm of Definition \ref{d:topomLF}, restricted to the space of molecules. It can be treated as a topometric analog of the 
 transportation cost norm. 
 It is also an attractive direction to study extreme points for the unit ball in $\mathcal{F}(\mathcal{M})$.

\sk 

Let $V$ be a normed space. 
Let us say that a normed subspace $W \subseteq V^*$ is \textit{norming} for $V$ if $||v||=\sup_{w\in W} \langle w,v \rangle$. The dual $V^*$ is always norming for $V$. By Definition \ref{d:topomLF}, $\Lip_0(\M)$ (which is $\F(\M)^*_{\t}$ by Theorem \ref{t:topomDual}) is norming for $\F(\M)$. We say that a Lipschitz function $f \colon (M,d) \to V$ is $W$-\textit{admissible} if $f \colon M \to (V, \s(V,W))$ is continuous (equivalently, $v^* \circ f \colon M \to \R$ is $\t$-continuous for every $v^* \in W$). Every weakly continuous function $f \colon M \to V$ is $W$-admissible for every $W \subset V^*$. 
The natural map $\d \colon M \to \F(\M)$ from Theorem \ref{t:topomDual} is $\Lip_0(M,d,\t)$-admissible.

\begin{theorem}[Universal property for topometric Lipschitz-free spaces] \label{t:univ} 
	
	Let $\M:=(M,d,\t,\0)$ be a pointed completely regular topometric space. Suppose that $(V,||\cdot||_V)$ is a Banach space and 
	$W \subseteq V^*$ is a norming subspace for $V$. Assume that 
	 $f \colon M \to V$ is a $W$-admissible function. 
%
	Then there exists a unique linear map $L_f \colon \mathcal{F}(\M) \to V$  such that $f=L_f \circ \d$ and $||L_f||=||f||_{\mathrm{Lip}}$.	
	\begin{equation*}
		\xymatrix { \mathcal{M} \ar[dr]_{f} \ar[r]^{\d} & \mathcal{F}(\mathcal{M})
			\ar[d]^{L_f} \\
			&  V \ \ & }
	\end{equation*} 
\end{theorem}
\begin{proof} Recall that $\{\d(x): x \in M \setminus \{\0\}\}$ is linearly independent in $\mathcal{F}(\mathcal{M})$.
	Thus, one may (uniquely) extend $f \colon M \to V$ to a linear function $L^m_f \colon \Mol(M) \to V$ such that $f=L^m_f \circ \d$, where $\Mol(M)=span (\d(M))$. Since $\Mol(M)$ is norm dense in $\mathcal{F}(M)$, it is enough to show $||f||_{\Lip}=||L^m_f||$. 
	Then we will have a desired canonical continuous linear extension on the completion $L_f \colon \F(\M) \to V$. 
	
	Clearly, $||f||_{\Lip} \leq ||L^m_f||$ (because $\d \colon M \to \Mol(M)$ is an isometric embedding). Therefore, we have only to veify that $||L^m_f|| \leq ||f||_{\Lip}$ on the space of molecules. 
	%
	%
	%
	
	Recall that by our definition, 
	$v^* \circ f$ is $\t$-continuous for every $v^* \in W$. 	
Thus, the restriction $v^* \circ L^m_f \circ \d$ of $v^* \circ L^m_f \colon \Mol(M) \to \R$ on $M$ is $\t$-continuous. 
Denote by $\psi_{v^*} \colon \F(\M) \to \R$ 
the canonical continuous linear extension of $v^* \circ L^m_f$. Then $\psi_{v^*} \circ \d$ is also continuous on $(M,\tau)$. Hence, $\psi_{v^*} \in   \mathcal{F}(\mathcal{M})_{\t}^{*}$. Moreover, the operator norms are the same. That is, $||\psi_{v^*}||=||v^* \circ L^m_f||$.  
	
	By Theorem \ref{t:topomDual} and (\ref{eq:form}),   $\mathcal{F}(\mathcal{M})_{\t}^{*} \to \Lip_0 (M,d,\t), \psi \mapsto \psi \circ \d$ is an isometric isomorphism. 
	This implies that $||v^* \circ L^m_f||=||v^* \circ L^m_f \circ \d||_{\Lip}$.  
	Using these observations (and the assumption that $W$ is norming for $V$), we get 
	$$
	\|L^m_f\|= \sup_{u \in B_{\Mol(M)}} ||L^m_f(u)||_V= \sup_{v^* \in B_{W}} \sup_{u \in B_{\Mol(M)}} |(v^*\circ L^m_f)(u)|=\sup_{v^* \in B_{W}} \|v^*\circ L^m_f\|=
	$$
	$$= 
	\sup_{v^* \in B_{W}}  \|v^* \circ L^m_f \circ \d\|_{\Lip} 
	=\sup_{v^* \in B_{W}} ||v^* \circ f||_{\Lip} \leq \|f\|_{\Lip}. 
	$$
\end{proof}

\begin{theorem}[Canonical linearization] \label{t:CanLinTopom} 
	
 For any Lipschitz-continuous (topometric morphism) map $f \colon \M_1:=(M_1,d,\t_1, \0) \to \M_2:=(M_2,d_2,\t_2, \0)$ between pointed topometric spaces, there exists a unique continuous linear map $L_f  \colon \mathcal{F}(\M_1) \to \mathcal{F}(\M_2)$ such that $L_f \circ \delta_1 = \delta_2 \circ f$ and $\|L_f\| = \|f\|_{\Lip}$.	
\end{theorem}
\begin{proof}
The natural map $\d \colon M \to \F(\M)$ from Theorem \ref{t:topomDual} is $\Lip_0(M,d,\t)$-admissible. Now, apply Theorem \ref{t:univ}  to the $\Lip_0(\M_2)$-admissible map $\delta_2 \circ f \colon M_1 \to \F(\M_2)$ and $V:=\F(\M_2)$. 	
\end{proof}

\begin{example} \label{ex:10metric}
	Let $(M,\tau)$ be a compact metrizable space with a distinguished point $0\in M$.
	Via the weak$^\ast$ embedding $x\mapsto\delta_x\in C(M)^\ast$, define
	\[
	d_\ast(x,y):=\|\delta_x-\delta_y\|=
	\begin{cases}
		0,& x=y,\\
		1,& \{x,y\}=\{0,z\}\ \text{with }z\neq0,\\
		2,& x\neq y\neq 0.
	\end{cases}
	\]
	Then for the compact topometric space $\mathcal M:=(M,d_\ast,\tau)$ we have:
	\begin{enumerate}
		\item $\displaystyle \Lip_0(M,d_\ast,\tau)=\{\,f\in C(M):f(0)=0\,\}$ and
		$\|f\|_{\Lip}=\|f\|_\infty$.
		In particular, $\Lip_0(M,d_\ast,\tau)$ is a closed subspace of $C(M)$ of codimension~1,
		hence a separable Banach space.
		\item $\mathcal F(\mathcal M)$ is isometrically isomorphic to $\ell^1(M\setminus\{0\})$.
		Consequently, $\mathcal F(\mathcal M)$ is separable if and only if $M$ is countable.
	\end{enumerate}
	This (with uncountable $M$) yields many examples with proper inclusion
	$\Lip_0(M,d,\tau)\subseteq \mathcal F(\mathcal M)^\ast$.
\end{example}

\begin{proof}
	(1) If $f\in B_{\Lip_0(M,d_\ast,\tau)}$ then for all $x,y\in M$,
	$|f(x)-f(y)|\le d_\ast(x,y)$, so $\|f\|_\infty\le1$ and $f(0)=0$.
	Thus the unit ball is $\{f\in C(M):f(0)=0,\ \|f\|_\infty\le1\}$.
	Conversely, for $f\in C(M)$ with $f(0)=0$,
	\[
	\|f\|_{\Lip}
	=\sup_{x\ne y}\frac{|f(x)-f(y)|}{d_\ast(x,y)}
	=\max\Bigl\{\sup_{x\ne0}\tfrac{|f(x)|}{1},\ \sup_{\substack{x\ne y\\ x,y\ne0}}\tfrac{|f(x)-f(y)|}{2}\Bigr\}
	=\|f\|_\infty.
	\]
	Hence $\Lip_0(M,d_\ast,\tau)=\{f\in C(M):f(0)=0\}$ with the supremum norm.
	
	(2) For $x\in M\setminus\{0\}$ set $e_x:=\delta_x-\delta_0\in \mathcal F(\mathcal M)$.
	Then $\|e_x\|=d_\ast(x,0)=1$.
	For any finitely supported $\{c_x\}\subset\mathbb R$,
	\[
	\Bigl\|\sum_x c_x e_x\Bigr\|
	=\sup\Bigl\{\Bigl|\sum_x c_x f(x)\Bigr|:\ f\in C(M),\ f(0)=0,\ \|f\|_\infty\le1\Bigr\}
	\le \sum_x |c_x|. 
	\]
	Equality holds by choosing $f$ with $f(x)=\operatorname{sgn}(c_x)$ on the finite support
	and $f(0)=0$ (Tietze extension), whence
	$\bigl\|\sum_x c_x e_x\bigr\|=\sum_x|c_x|$.
	Therefore $\{e_x\}_{x\in M\setminus\{0\}}$ is an isometric $\ell^1$-basis of
	$\mathcal F(\mathcal M)$, so $\mathcal F(\mathcal M)\cong \ell^1(M\setminus\{0\})$.
\end{proof}

\sk 
\begin{example} \label{e:ForTopometricRevised}  Consider the $\tau$-compact space 
(see Example \ref{ex:NotBounded}) 
$( M = \{0\} \cup \left\{ \frac{1}{n} : n \in \mathbb{N} \right\}$ with the metric $d \left( \tfrac{1}{n}, \tfrac{1}{m} \right) = |n - m|, \ d\left( \tfrac{1}{n}, 0 \right) = n.$ 
For this topometric space $\mathcal{M} = (M, d, \tau)$  
we have: 
		\begin{enumerate}
			\item $\Lip_0(\mathcal{M})$ is separable and not complete, where $$\Lip_0(\mathcal{M}) = \left\{ a = (a_n)_{n \in \mathbb{N}_0} : a_0 = 0, \; \lim_{n \to \infty} a_n = 0 \right\}, 
			\ \|a\| := \sup_n |a_{n+1} - a_n|;$$  			
			\item the linear isometry $\mathcal{F}(\mathcal{M}) \cong \ell^1$. 
		\end{enumerate}
			\end{example}
		\begin{proof} (1) 
		Writing \( a_n = f\left( \tfrac{1}{n} \right) \) for \( n \in \mathbb{N} \) and \( a_0 = f(0) \), every function \( f \in \Lip_0(\mathcal{M}) = \Lip_0(M,d,\tau)\) corresponds to a real sequence \( (a_n)_{n \in \mathbb{N}_0} \) with \( a_0 = 0 \) and the property
		\[
		\|f\|_{\Lip} = \sup_{n \in \mathbb{N}_0} |a_{n+1} - a_n| < \infty.
		\]
		This can be verified by the telescopic sum argument  (since $a_n-a_0=\sum_{j=0}^{n-1}(a_{j+1}-a_j)$),  
		 leading to the identification 
		\[
		\Lip_0(\mathcal{M}) = \left\{ a = (a_n)_{n \in \mathbb{N}_0} : a_0 = 0, \; \lim_{n \to \infty} a_n = 0 \right\},
		\quad \|a\| := \sup_n |a_{n+1} - a_n|.
		\]
		
		\smallskip 
Note that, norm convergence $\|a^{(k)}-a\|\to0$ implies coordinatewise convergence: for each fixed $n$,
		$|a^{(k)}_n-a_n|\le n\,\|a^{(k)}-a\|\to 0$.
		
		\sk 
		This normed space \( V := \Lip_0(\mathcal{M}) \) is \textbf{not complete}. For instance, consider the following sequence 
		of vectors \( (a^{(k)})_{k \in \mathbb{N}} \subset V \) as follows. For each \( k \in \mathbb{N} \), define the sequence \( a^{(k)} = (a_n^{(k)})_{n \in \mathbb{N}_0} \) by:
		\[
		a^{(k)}_0 = 0, \quad
		a^{(k)}_n =
		\begin{cases}
			1, & \text{if } 1 \le n \le k, \\
			\frac{2k - n}{k}, & \text{if } k < n \le 2k, \\
			0, & \text{if } n > 2k.
		\end{cases}
		\] 
		Then this sequence is Cauchy in $V$. Indeed, 
		a direct computation, for 
		$(\Delta a)_n:=a_{n+1}-a_n \ (n\in\mathbb{N}_0)$  gives, for every $k$,
		$$
		(\Delta a^{(k)})_0=1,\qquad
		(\Delta a^{(k)})_n=
		\begin{cases}
			0,& 1\le n\le k-1,\\[1mm]
			-\dfrac{1}{k},& k\le n\le 2k-1,\\[2mm]
			0,& n\ge 2k.
		\end{cases}
		$$
		Hence, for $k,m\in\mathbb{N}$,
		$$
		\|a^{(k)}-a^{(m)}\|
		=\|\Delta a^{(k)}-\Delta a^{(m)}\|_\infty
		=\max\Bigl\{\frac{1}{k},\frac{1}{m}\Bigr\}\xrightarrow[\min\{k,m\}\to\infty]{}0,
		$$ 
		so $\bigl(a^{(k)}\bigr)$ is Cauchy. 
		
		However, it is not norm convergent in $V$.
		Assume towards a contradiction that \(a^{(k)} \to a\) in norm. Then, in particular,
		\(a^{(k)}_n \to a_n\) for every \(n\) (coordinatewise).
		By our computation \(\|\Delta a^{(k)}-e_0\|_{\infty}=\frac{1}{k} \to 0\).
		Since \(\Delta:V\to\ell_\infty\), \((\Delta a)_n=a_{n+1}-a_n\), is an isometry,
		we also have \(\Delta a^{(k)}\to \Delta a\) in \(\ell_\infty\), hence \(\Delta a=e_0\),  
		where \(e_0=(1,0,0,\ldots)\).
		Therefore, for all \(n\ge 0\),
		\[
		a_{n+1}-a_n=(\Delta a)_n=(e_0)_n=
		\begin{cases}
			1,& n=0,\\
			0,& n\ge 1.
		\end{cases}
		\]
With \(a_0=0\) this yields \(a=(0,1,1,1,\ldots)\). 	So, \(a^{(k)}\) converges \textit{pointwise} to a vector  $a:=(0,1,1,1,\cdots)$ which lies outside \( \Lip_0(\mathcal{M}) \) (see Remark \ref{r:Kizmaz}). Thus \(\bigl(a^{(k)}\bigr)\) has no norm limit in \(V\).

\sk 
\textbf{Separability} of $V$. 
The map $\Delta:V\to\ell_\infty$ 
is an isometry onto
$H:=\{g\in\ell_\infty:\sum_{j=0}^{n-1} g_j\to0\}$.
For each $g\in H$ and $\varepsilon>0$, one can approximate $g$ uniformly by rational, finitely supported, zero-sum sequences.  
Hence $H$ is separable, and so are $V$ and its completion.  

\sk 
	
		(2) Recall that $\mathcal{F}(\mathcal{M})$ is the closed span of the molecules 
		$\delta_n$ with $\delta_0=\0$. 
		Every finitely supported molecule has the form
		\[
		\mu=\sum_{i=1}^N c_i \delta_{n_i}, \qquad \sum_{i=1}^N c_i=0.
		\]
		By Definition~\ref{d:topomLF}, the norm is
		\[
		\|\mu\|=\sup\Bigl\{\sum_{i=1}^N c_i f(n_i) : f\in \Lip_0(M,d,\tau),\ \|f\|_{\Lip}\le 1\Bigr\}.
		\]
		
		Since $\Lip_0(M,d,\tau)$ is not complete, it is a proper subspace of $\Lip_0(M,d)$. 
		However, for a finitely supported $\mu$ only the values $f(0),f(1),\dots,f(K)$ with 
		$K=\max_i n_i$ are relevant in the supremum. Beyond $K$ one may redefine $f$ so that 
		$f(n)\to f(0)$ as $n\to\infty$, without changing the values $f(n_i)$, and hence without 
		changing the pairing with $\mu$. 

So, we have the linear isometry $\mathcal{F}(\mathcal{M})=\mathcal{F}(M,d,\tau) \cong \mathcal{F}(M,d)=\mathcal{F}(\mathbb{N}_0)$. On the other hand, 
it is a classical fact (see, e.g., Aliaga \cite[Example~2.3]{Aliaga20}) that $\mathcal{F}(\mathbb{N}_0) \;\cong\; \ell^1.$  
	\end{proof}	

\sk   
\begin{remark}[Relation to Kizmaz space] \label{r:Kizmaz} 
	In fact, the sequence \(\bigl(a^{(k)}\bigr)\) from \(V\), discussed in the proof of
	Example~\ref{e:ForTopometricRevised}, converges to the vector \(a\) in the Kizmaz
	space \(\ell_{\infty}(\Delta)\) (see~\cite{Kiz}), where
	\[
	\ell_{\infty}(\Delta):=\bigl\{\,b=(b_n)_{n\in\mathbb N}:\ \sup_n|b_{n+1}-b_n|<\infty\,\bigr\},
	\qquad
	\|b\|_{\Delta}:=|b_1|+\sup_n|b_{n+1}-b_n|.
	\]
	Define \(T:V\to \ell_{\infty}(\Delta)\) by
	\(T(x_0,x_1,x_2,\dots)=(x_1,x_2,x_3,\dots)\). Then
	\[
	\|a\|_{V}\ \le\ \|T(a)\|_{\Delta}\ \le\ 2\,\|a\|_{V}
	\quad(\text{since }|a_1|\le\|a\|_{V}\text{ and }\sup|\Delta T(a)|=\sup|\Delta a|),
	\]
	so, \(T\) is a linear topological embedding of \(V\) into the Banach space
	\(\ell_{\infty}(\Delta)\) (as sets one has \(T(V)=c_0\)). In particular, the completion
	\(\widehat{V}\) embeds linearly and topologically into \(\ell_{\infty}(\Delta)\).
\end{remark}

\begin{proposition} \label{p:compBound}
If $\M:=(M,d,\t,\0)$ is a pointed topometric space with compact $\t$ and bounded $d$, then the normed space  $\Lip_0(M,d,\t)$ is complete.  	
\end{proposition} 
\begin{proof}
	Let \( V := \Lip_0(M,d,\tau) \), where \( (M,d,\tau,\0) \) is a pointed topometric space with compact topology \( \tau \) and bounded metric \( d \leq D \). Take a Cauchy sequence \( (f_n) \subset V \) in the Lipschitz norm \( \|\cdot\|_{\Lip} \). 
	
	Since \( \|f_n - f_m\|_\infty \le D \cdot \|f_n - f_m\|_{\Lip} \), the sequence \( (f_n) \) is also Cauchy in \( (C(M), \|\cdot\|_\infty) \), hence converges uniformly to some \( f \in C(M) \). Moreover, each \( f_n \in \Lip_0(M,d) \), and the pointwise limit of a uniformly convergent sequence of Lipschitz functions with uniformly bounded Lipschitz constants is again Lipschitz with the same bound. So \( f \in \Lip_0(M,d) \). 
	
	Since each \( f_n \in C(M,\tau) \) and the limit is uniform, it follows that \( f \in C(M,\tau) \) and \( f(0) = 0 \). Hence \( f \in \Lip_0(M,d,\tau) \). Finally, \( \|f_n - f\|_{\Lip} \to 0 \), so \( f_n \to f \) in \( V \), proving completeness.
\end{proof}

\begin{remark} \label{r:BALLS} 
It would be interesting to study topometric Lipschitz-free spaces coming from Banach spaces. For example, what can be said about $\F(\M)$ and $\Lip_0(M,d,\tau)$ for the following two topometric spaces canonically defined for every Banach space $V$:  (a) $\M:=(B_V,d_{norm},\tau_{weak})$; (b) $\M:=(B_{V^*},d_{norm},\tau_{weak^*})$. In case (a), $\Lip_0(M,d,\tau)$ naturally contains $V^*$ as a subspace and in case (b), it contains $V$.  Precise description seems to be hard. A particular but important case is the case of $V=\ell_2$. 
 Proposition \ref{p:compBound} guarantees that $\Lip_0(B_{V^*},d_{norm},\tau_{weak^*})$ is always complete. 
Note that $\Lip_0(B_V,d_{norm},\tau_{weak})$ might be nonseparable for separable $V$. Indeed, this happens for $V:=\ell_1$ (for every separable $V$ which is not Asplund). 	
\end{remark}

\sk  
\section{Induced linear isometric group actions} 
\label{s:actions} 

We begin by recalling some standard notions concerning group actions and $G$-compactifications.  

By a $G$-\textit{space}, we mean a topological space $X$ endowed with a continuous action $\pi \colon G \times X \to X$, where $\pi(g,x) = gx$ 
of a topological group $G$ on $X$. 
If $G$ is a discrete group then the action is continuous if and only if every $g$-translation $t^g \colon X \to X, t^g(x)=gx$ is continuous; if and only if $\pi \colon G_{\mathrm{disc}} \times X \to X$ is continuous, where $G_{\mathrm{disc}}$ is the discrete copy of $G$.  

If $X$ is a  compact $G$-space, sometimes we say that $X$ is a \textit{$G$-flow}. If $X$ is a convex compact set in a locally convex space and the continuous action of $G$ is by affine maps, then we say that $X$ is an \textit{affine $G$-flow}. 
A continuous map $f \colon X_1 \to X_2$ between two $G$-spaces is called a $G$-\textit{map} (or \textit{equivariant}) if $f(gx) = g f(x)$ for all $g \in G$ and $x \in X_1$.  

\subsection*{$G$-compactifications} 

We emphasize that in the present paper ``compactification" is not necessarily a topological embedding.  
More precisely, a continuous dense map $\nu \colon X \to Y$ into a compact Hausdorff space $Y$ is called a \textit{compactification} of $X$. If $X$ and $Y$ are $G$-spaces and $\nu$ is a $G$-map, then $\nu$ is said to be a $G$-\textit{compactification}. If $\nu$ is also a topological embedding, it is called \textit{proper}. Recall that $\nu$ is proper iff the canonically induced Banach unital subalgebra $\mathcal{B}_{\nu}:=\{f \circ \nu: f \in C(Y)\}$ of $C_b(X)$ separates the points and closed subsets. 

Recall that two compactifications $\nu_1 \colon X \to Y, \nu_2 \colon X \to Y_2$ are equivalent if there exists a homeomorphism $h \colon Y_1 \to Y_2$ such that $\nu_2=h \circ \nu_1$. Or, iff these compactifications have the same induced Banach algebras $\mathcal{B}_{\nu_1}=\mathcal{B}_{\nu_2}$. If $X$ is a $G$-space and $\nu_1, \nu_2$ are equivalent $G$-compactifications then the homeomorphism $h$ above necessarily is a $G$-map.  
The Samuel compactification of a uniform space $(X,\mathcal{U})$ is the compactification induced by the algebra $\Unif^b(M,\mathcal{U})$ of all $\mathcal{U}$-uniformly continuous bounded real functions on $X$.  
The following technical result probably is well known. However we could not find a direct reference. 
\begin{fact} \label{f:SW} 
	Let $F \subseteq C_b(X)$ be a set of continuous bounded functions on $X$. Denote by 
	$$
	\nu_F \colon X \to \R^F, \ \ x \mapsto (f(x))_{f \in F}
	$$ 
	the diagonal function and by $Y:=cl_p(\nu_F(X))$ (necessarily compact) subset of $\R^F$.  
	Then the Banach algebra $\{f \circ \nu: f \in C(Y)\}$ associated to the induced compactification $\nu_F \colon X \to Y$ is the smallest unital Banach subalgebra $\mathcal{A}_F$ of $C_b(X)$ which contains $F$.  
\end{fact}
\begin{proof}
	For the compactification $\nu_F \colon X \to Y$ we have the induced inclusion of algebras $\nu_F^* \colon C(Y) \hookrightarrow C_b(X)$. 
	Then $\nu_F^*(p_f)=f$, for every $f \in F$, where $p_f \colon Y \to \R$ is the corresponding coordinate projection and it extends $f \colon X \to \R$. Thus, $\mathcal{A}_F \subseteq \nu_F^*(C(Y))$. On the other hand,   
	the family of all projections $P:=\{p_f \colon Y \to \R : f \in F\}$ separate the points of the compact space $Y$. Therefore, by the Stone-Weierstrass theorem the unital subalgebra generated by the subset $P$ is just $C(Y)$. So, $\mathcal{A}_F \supseteq  \nu_F^*(C(Y))$ and we conclude that $\mathcal{A}_F = \nu_F^*(C(Y))$.     	
\end{proof}

\sk 
An action of $G$ on a metric space $(M,d)$ is said to be \textit{isometric} if every translation map $t^g \colon M \to M$, defined by $t^g(x) := gx$, is an isometry. 

We use below several times the following simple technical result.  
\begin{fact} \label{f:cont} 
	An isometric action $\pi \colon G \times M \to M$ is continuous if and only if the orbit map $\operatorname{orb}_y \colon G \to M$, $g \mapsto gy$, is continuous for every $y\in Y$ in some dense subset $Y \subset M$. 	
\end{fact}

As before, let $(M,d,\0)$ denote a pointed metric space. Suppose $\pi \colon G \times M \to M$ is a continuous isometric action such that $g\0=\0$ for all $g \in G$. Recall the canonical isometric embedding  
\[
\delta \colon M \hookrightarrow \mathrm{Mol}(M,d), \quad x \mapsto \delta_x.
\]

By linearly extending the original action $\pi$ from the set $\delta(M)$ to the normed space $\mathrm{Mol}(M)$ of all molecules, we obtain a linear action  
$G \times \mathrm{Mol}(M) \to \mathrm{Mol}(M)$.  
It is straightforward to verify (by the universal property Fact \ref{f:facts}.4) that this action is isometric and  separately continuous, and hence jointly continuous by Fact~\ref{f:cont}. Moreover, passing to the completion, we obtain a uniquely determined continuous linear isometric action
\[
G \times \mathcal{F}(M) \to \mathcal{F}(M). 
\]

\begin{remark} \label{r:IsoFix} 
If $G \times M \to M$ is an isometric action with $(M,d)$ not necessarily pointed (and not necessarily containing a $G$-fixed point), then one may try to adjoin a new point $\0$ which will be $G$-fixed and the extended action of $G$ on $M^+:=M \cup \{0\}$ will remain isometric. It is easy if $(M,d)$ is bounded. Indeed, we can define $d^+(\0,x)=c_0 \geq \operatorname{diam}(M,d)$. This fact is well known and easy to verify. See, for example, \cite[Proposition 2.10]{Me-cs07}. 
Moreover, an exact criteria was obtained by Schr\"{o}der \cite{Sc}. 
It asserts that a metric space $(M,d)$ can be extended by adding a $G$-fixed point getting again an isometric action if and only if 
all orbits (equivalently, if one of the orbits) 
$Gx$ are bounded for every $x \in M$. In fact, all this is true for extension of monoid $1$-Lipschitz (i.e., non-expansive) actions.  
\end{remark}

Let $h_0 \colon G \to \mathrm{Is}_{\text{lin}} (\mathcal{F}(M))$ be the  canonically defined continuous group homomorphism, where $\mathrm{Is}_{\text{lin}} (\mathcal{F}(M))$ is the topological group of all \textbf{linear} isometries endowed with the \textit{strong operator topology} (SOT). This is the  topology inherited from the product 
$(\mathcal{F}(M), norm)^{\mathcal{F}(M)}$. Similarly, the topology on $\mathrm{Is}_{\text{lin}} (\mathcal{F}(M))$ inherited from 
$(\mathcal{F}(M), weak)^{\mathcal{F}(M)}$ is said to be the \textit{weak operator topology} (WOT).

Proposition \ref{t:light} below implies that SOT and WOT coincide on the subgroup $h(G) \subset \mathrm{Is}_{\text{lin}} (\mathcal{F}(M))$. 

Denote by $\Is(M)$ the group of all self-isometries onto of $M$ (fixing $\0$, whenever $M$ is pointed). It is a topological group with the pointwise topology. 

\begin{lemma} \label{l:emb} 
$h \colon \Is(M) \hookrightarrow (\mathrm{Is}_{\text{lin}} (\mathcal{F}(M)), SOT)$ is an embedding of topological groups. 	
\end{lemma}
\begin{proof}
Indeed, as we already explained, 
the linear action of $G=\Is(M)$ on $\mathcal{F}(M)$ is continuous. Therefore, $h$ is continuous (where $\mathrm{Is}_{\text{lin}} (\mathcal{F}(M))$ carries SOT). 
Moreover, the restricted continuous action of the image $h(\Is(M))$ on $\d(M) \subset \mathcal{F}(M)$ is an equivariant copy of the original action $\pi$. Thus, $h$ is injective and 
every orbit map $orb_{\d(x)} \colon h(\Is(M)) \to \d(M)$ is norm continuous for every $x \in M$. Hence $h$, in fact, is an embedding of topological groups.  
\end{proof}
 



	Let $V$ be a Banach space. Recall that a subgroup $G$ of $\mathrm{Is}_{\text{lin}} (V)$ is said to be \emph{light} (see \cite{Me-OperTop01, Me-FRepres03}) if $SOT$ and $WOT$ agree on $G$. 

\begin{proposition} \label{t:light} 	
	Let 
	$G \times M \to M$ be an isometric continuous action of a topological group $G$ on a pointed metric space $M$. Then the weak continuity of a homomorphism $h \colon G \to \mathrm{Is}(M) \subset \mathrm{Is}_{\text{lin}} (\mathcal{F}(M))$ implies its strong continuity. In particular, WOT and SOT on $\mathrm{Is}(M)$ coincide. That is, the subgroup $\mathrm{Is}(M,d)$ is light in $\mathrm{Is}_{\text{lin}} (\mathcal{F}(M))$.
\end{proposition}
\begin{proof} 
	Let $h \colon G \to \mathrm{Is}(M) \subset \mathrm{Is}_{\text{lin}} (\mathcal{F}(M))$ be weakly continuous. That is, the orbit map $orb_v \colon G \to \mathcal{F}(M)$ is weakly continuous for every $v \in \mathcal{F}(M)$. 
	By Proposition \ref{p:w=n}.1,  
	weak and norm topologies coincide on $M =\delta(M) \subset \mathcal{F}(M)$. Hence, for every $v \in M$ the orbit maps 
	$orb_v \colon G \to \mathcal{F}(M)$ are norm continuous (because $Gv \subseteq M$). By the continuity of linear operations in $\mathcal{F}(M)$ it is also clear that $orb_v$ are norm continuous for every vector $v=\sum_{i=1}^n c_i \d_{m_i}$  from the linear span of $M$. 
	That is, for every $v \in \Mol(M)$.  
	Since $\Mol(M)$ is norm dense in $\mathcal{F}(M)$, by 
	Fact \ref{f:cont}, we obtain that 
	the action $G \times \mathcal{F}(M) \to \mathcal{F}(M)$ is continuous. 
	 This implies that $h$ is norm continuous.     
\end{proof}
 
 \begin{question} \label{q:light} 
 For which pointed metric spaces $M$ the group 	$\mathrm{Is}_{\text{lin}} (\mathcal{F}(M))$ is light ? 
 \end{question}

\begin{remark} \label{r:light}  
 Recall that 
	$\mathrm{Is}_{\text{lin}} (V)$ is light for every reflexive Banach space $V$ \cite{Me-OperTop01, Me-FRepres03}, while the group  $\Is_{\text{lin}}(C([0,1]^2))$ is not light. We refer to \cite{AFGR} for more information which contains also several additional examples and counterexamples. For example, $\Is_{\text{lin}}(L_1[0,1])$ is not light.  
	
	Since the Lipschitz-free space $\mathcal{F}(\R)$ is $L_1[0,1]$, it follows that $\mathrm{Is}_{\text{lin}} (\mathcal{F}(M))$ need not be light in general. Surprizingly enough, the class of pointed metric spaces $M$ with light $\mathrm{Is}_{\text{lin}} (\mathcal{F}(M))$ is quite large and contains the so-called \emph{weak Prague spaces $M$}; see  \cite[Proposition 6.2 and Remark 6.3]{CDT}. 
\end{remark}

\sk 
\subsection{Induced dual action} As before, let $\pi \colon G \times M \to M$ be a continuous action by isometries with fixed $\0$. We have the corresponding continuous isometric linear action 
$$G \times \mathcal{F}(M) \to \mathcal{F}(M).$$
It leads to the induced \textit{dual action} on the dual space  $\mathcal{F}(M)^*=\Lip_0(M)$ 
$$G \times \Lip_0(M) \to \Lip_0(M), \ \ (g \varphi)(v):=\varphi(g^{-1}v)$$  
by linear isometries. This dual action need not be norm continuous even for compact $G$ (Example \ref{ex:NotCont}).  
 Note that if a Banach space $V$ is Asplund then the dual action $G \times V^* \to V^*$ is continuous for every linear continuous isometric action of a topological group $G$ on $V$. 
  However, $\F(M)$ typically contains an isomorphic copy of $\ell_1$ (so, it is not Asplund). Compare Question \ref{q:SEPARABtopom}. 

According to the following result, the weak-star topology gives a rich and important source of continuous actions on 
any bounded $G$-invariant subset of $V^*$. 

\begin{fact} \label{f:DualActionGen} 
	Let $V$ be a normed space and $\pi \colon G \times V \to V$ is a continuous action by linear isometries. 
	Then  
	the induced dual action $\pi^* \colon G \times B_{V^*} \to B_{V^*}$ is continuous, where $B_{V^*}$ is the weak-star compact unit ball in the dual space $V^*$. Remains true for every bounded 
	$G$-invariant subset of $V^*$.   
\end{fact}
\begin{proof} 
	The reason is that the induced homomorphism $(\mathrm{Is}_{\text{lin}} (V), SOT) \to \Homeo(B_{V^*})$ is an embedding of topological  groups, where $\Homeo(B_{V^*})$ carries the compact-open topology.   
	See more details, for example, in \cite[Lemma 2.4]{Me-cs07}.  	
%
%
%
%
\end{proof}

\begin{corollary} \label{l:DualAction} 
	Let $\pi \colon G \times M \to M$ be a continuous action by isometries with fixed $\0$.
	The induced dual action $\pi^* \colon G \times B_{\mathcal{F}(M)^*} \to B_{\mathcal{F}(M)^*}$ is continuous, where $B_{\mathcal{F}(M)^*}$ is the weak-star compact unit ball in the dual space ${\mathcal{F}(M)^*}$. This remains true for every bounded 
	 $G$-invariant subset of $\mathcal{F}(M)^*$.   
\end{corollary}

\begin{definition} \label{d:affineComp} 
	(see \cite{GM-fixp12}) 
	An {\em affine \(G\)-compactification} of a \(G\)-space \(X\) is a
	pair \((\a,Q)\) where \(\a\colon X\to Q\) is a continuous \(G\)-map,
	\(Q\) a compact convex affine \(G\)-flow in a locally convex linear space, and
	\(\overline{\mathrm{co}}(\a(X))=Q\).
\end{definition}

 In terms of Fact \ref{f:DualActionGen} 
 let $f \in V^*$. Consider the induced $G$-compactification of the orbit $Gf$: 
$$
K_f:=cl_{w^*}(Gf)=\overline{{Gf}}^{w^*}.
$$  
and also the induced affine $G$-compactification: 
 $$
 Q_f:=cl_{w^*} co(K_f)=\overline{co}^{w^*}(Gf). 
 $$ 		
Then $K_f$ is a $G$-flow and $Q_f$ is an affine $G$-flow by Fact  \ref{f:DualActionGen}. 
For a related material see below Question \ref{q:ind}.2 and Remark \ref{r:amenability}.

\sk 
\subsection{Induced double dual action} 
The original continuous isometric action 
$\pi \colon G \times M \to M$ implies (not necessarily continuous) linear isometric action $G \times \Lip_0(M) \to \Lip_0(M)$, which in turn, induces the (dual) action by linear isometries 
$$
	G \times \Lip_0(M)^* \to \Lip_0(M)^* 
$$
on $\Lip_0(M)^*$ (which is the bidual $\mathcal{F}(M)^{**}$).  
Every $g$-translation 
$$t^g \colon (\Lip_0(M)^*,w^*) \to (\Lip_0(M)^*,w^*)$$
 is weak-star continuous. Therefore, $(\Lip_0(M)^*,w^*)$ is a $G_{\mathrm{disc}}$-space, where $G_{\mathrm{disc}}$ is the discrete copy of $G$.   
 In contrast to the dual action on $\Lip_0(M)$ (remember Corollary \ref{l:DualAction}),  
 for this case, the continuity of the restricted action on weak-star compact (bounded) subsets of the bidual is not guaranteed in general 
 if $G$ is not discrete. 
  
 Consider on the dual space $\Lip_0(M)^*$ the ``weak-star uniformity" $\mathcal{U}_*$. That is, the (weak) uniformity $\mathcal{U}_*$ generated by the collection $\{f \colon \Lip_0(M)^* \to \R: f\in \Lip_0(M)\}$. 
 The topology of $\mathcal{U}_*$ is just the usual weak-star topology $w^*$.

Recall that    
$\d \colon M \hookrightarrow (\Lip_0(M)^{*},weak^*)$ is a topological embedding (Proposition \ref{p:w=n}).  
\begin{remark} \label{r:APS} 
	A recent result \cite[Proposition 2.3]{APS24} shows that the so-called \emph{Lipschitz realcompactification} $M^{\mathcal{R}}$ \cite{GarMer} of $(M,d)$ can be identified topologically with the subset $(\overline{\d(M)}^{w^*}, w^*)$ of $\Lip_0(M)^{*}$. 
\end{remark}

 Every $t^g$-translation $\Lip_0(M)^* \to \Lip_0(M)^*$ is $\mathcal{U}_*$-uniformly continuous.  
Since $\d(M)$ is a $G$-invariant subset, its weak-star closure $$
M^{\mathcal{R}}:=(\overline{\d(M)}^{w^*},w^*) \subset (\mathcal{F}(M))^{**}
$$ is also $G$-invariant. The action of $G$ on $M^{\mathcal{R}}$ is at least $G_{\mathrm{disc}}$-continuous.

 \begin{question} \label{q:delta-starComp} 
	Study properties of the dense embedding (which is a $G_{\mathrm{disc}}$-compactification for bounded $d$)
	$\d_* \colon M \hookrightarrow M^{\mathcal{R}} \subset (\mathcal{F}(M))^{**}.$ 
	In particular, when the $G$-action on $M^{\mathcal{R}}$ is continuous ?  
\end{question}

In general, the answer to Question \ref{q:delta-starComp} is in the negative even for compact groups $G$ and bounded $d$ 
(Example \ref{ex:NotCont}). For a positive example, see Proposition \ref{p:positive}. 

If $(M,d)$ is bounded, then $\d_*$ is equivalent to the Samuel compactification of $(M,\Unif(d))$, where $\Unif(d)$ is the uniform structure of $d$.
 Indeed, for bounded $d$, $\Lip_0(M)$ is uniformly dense in the algebra $\Unif^d(M,d,\0)$ of all $d$-uniformly continuous bounded real functions which vanish at $\0$.

\sk 
In order to control $G$-compactifications we need to deal with a special class of RUC functions.

\sk 
\begin{definition} \label{d:uniform} 
	(See, for example, \cite{Vr-Embed77,Me-opit07})  
Let $\pi \colon G \times X \to X$ be a continuous action of a topological group $G$. 
\begin{enumerate}
	\item 
	A real continuous function $f \colon X \to \R$ is \emph{right uniformly continuous} 
	if the following holds:  
	$$
	\forall \eps>0 \ \exists U(e) : \ \ |f(ux)-f(x)|\leq \eps \ \ \  \text{for all} \ x \in X \ \text{and} \  u \in U(e), 
	$$ 	
	where $U(e)$ is a neighbourhood of the neutral element $e$ in $G$. Notation: $f \in \mathrm{RUC}_G(X)$. The subfamily of all \emph{bounded} RUC, we denote by $\mathrm{RUC}^b_G(X)$. 
	\item Let $(X,\mathcal{U})$ be a uniform space. We say that the action $\pi$ is \textbf{equiuniform} if every  $g$-translation $t^g \colon X \to X$ is  $\mathcal{U}$-uniformly continuous and 
	$$
	\forall \eps \in \mathcal{U} \  \exists U(e): 
	\ \ (ux,x) \in \eps \ \ \  \text{for all} \ x \in X \ \text{and} \ u \in U(e).
	$$ 
	
	This definition of equiuniformity is equivalent to the following requirement  
	$$
	\forall \eps \in \mathcal{U} \  \exists U(g_0): 
	\ \ (gx,g_0x) \in \eps \ \ \  \text{for all} \ x \in X \ \text{and} \ g \in U(g_0), 
	$$ 
	where $U(g_0)$ is a neighbourhood of the element $g_0 \in G$. 	
\end{enumerate} 
\end{definition}

Definition \ref{d:uniform}.2 appears in \cite{Br} and \cite{Vr-Embed77} under the names: \emph{motion equicontinuous} and \emph{``bounded uniformity}". 

\begin{remark} \label{r:equivalence} 
It is well known (see for example \cite{Vr-Embed77}) that  $\mathrm{RUC}^b_G(X)$ is a unital Banach subalgebra of $C_b(X)$ and the corresponding Gelfand compactification $\beta_G \colon X \to \beta_G X$ is the \emph{greatest $G$-compactification} ($G$-analog of Stone-\v{C}ech compactification) of $X$. Moreover, there exists a natural 1-1 correspondence between unital $G$-invariant closed subalgebras of $\mathrm{RUC}^b_G(X)$ and $G$-compactifications of $X$.  

Recall that the greatest $G$-compactification $\beta_G \colon X \to \beta_G X$ is not necessarily proper (even for Polish $G$ and $X$) \cite{Me-Ex88}. $\beta_G X$ might be even a singleton for nontrivial $X$ (Pestov \cite{Pest-Smirnov}). 
\end{remark}

\begin{fact} \label{f:G-compactFacts} \cite[Lemma 4.5]{Me-cs07} 
\begin{enumerate}
	\item 
	Let $(Y,\mathcal{U})$ be a uniform space 
and let $\pi \colon G \times Y \to Y$ be an action with uniform $g$-translations. Suppose that there exists a $G$-invariant dense subset $X \subseteq Y$ such that the inherited action $G \times X \to X$ is 
$\mathcal{U}|_X$-equiuniform. Then the original action $\pi$ on $Y$ is $\mathcal{U}$-equiuniform  and continuous. 
	\item Let $\pi \colon G \times X \to X$ be a continuous $\mathcal{U}$-equiuniform action on a uniform space $(X,\mathcal{U})$. Then the canonically extended 
$G$-action on the completion
 $$\widehat{\pi} \colon G \times \widehat{X} \to \widehat{X}$$ 
 is $\widehat{\mathcal{U}}$-equiuniform and continuous.  
%
\end{enumerate}	
\end{fact}
\begin{proof}  
(1)
We show that the action $\pi \colon G \times Y \to Y$
is $\mathcal{U}$-equiuniform (the continuity of $\pi$ is an easy corollary).  
Let $g_0 \in G$ and $\eps \in {\mathcal{U}}$. There exists
an entourage $\eps_1 \in {\mathcal{U}}$ such that $\eps_1 \subset \eps$ and
$\eps_1$ is a \textit{closed} subset of $Y \times Y$. 
Since $\pi_X \colon G \times X \to X$ is $\mathcal{U}|_X$-equiuniform,   
one may choose a neighborhood 
$U(g_0)$ of $g_0$ such that $(g_0x,gx) \in \eps_1$ for every $g \in U(g_0)$ and $x \in X$. For a given $y \in Y$ 
choose any net $(x_i)$ in $X$ which tends to $y$ in $Y$. Then for any given pair $(g_0, g)$ we have $\lim_i g_0x_i=g_0y$ and $\lim_i gx_i=gy$. Since $\eps_1$ is closed, we obtain that $(g_0y,gy) \in \eps$ for every $g \in U(g_0)$ and $y \in Y$. 
	
(2) Directly follows from (1) with $Y:=\widehat{X}$. 
\end{proof}

\begin{theorem} \label{t:positive} 
	Let $(M,\rho)$ be a pointed 
	metric space with a continuous isometric action of a topological group $G$. 
	\begin{enumerate}
		\item 	Assume that $\Lip_0(M) \subseteq \mathrm{RUC}_G(M)$. 
		Then the natural 
		action   
		$G \times M^{\mathcal{R}} \to M^{\mathcal{R}}$  
		on the Lipschitz realcompactification $ M^{\mathcal{R}}$ 
		is $G$-continuous.  
		\item If, in addition, $\rho$ is bounded, then 
		$\d_* \colon M \hookrightarrow M^{\mathcal{R}}$ 
		is a $G$-compactification if and only if $\Lip_0(M) \subseteq \mathrm{RUC}_G(M)$.  
	\end{enumerate} 
\end{theorem}
\begin{proof} 
	(1) 
	By Remark \ref{r:APS} we have an identification  $M^{\mathcal{R}}=(\overline{\d(M)}^{w^*}, w^*)$.  
	As we already mentioned above,  
	$\Lip_0(M)$ generates the weak-star uniformity $\mathcal{U}_*$ on $\Lip_0(M)^*$ and 
	$G \times M^{\mathcal{R}} \to M^{\mathcal{R}}$ is $G_{\mathrm{disc}}$-continuous. That is, the $g$-translations $M^{\mathcal{R}} \to M^{\mathcal{R}}$ are continuous. Since $\Lip_0(M) \subseteq \mathrm{RUC}_G(M)$, 
	the induced action on subspace uniformity $\mathcal{U}_*|_M$ is equiuniform. 
	Also, we know that $M$ is a dense $G$-invariant subspace of $M^{\mathcal{R}}$. Therefore, one may apply Fact \ref{f:G-compactFacts}.1.       
	
	(2) $M^{\mathcal{R}}$ is compact if (and only if) $(M,\rho)$ is bounded. One direction follows from assertion (1) because  $\Lip_0(M) \subseteq \mathrm{RUC}_G(M)$ then the action is continuous and hence we have a $G$-compactification. For the second direction use Remark \ref{r:equivalence} which implies that if the $G$-action on the compact space $M^{\mathcal{R}}$ is continuous then every $f \in \Lip_0(M)$ belongs to $\mathrm{RUC}_G(M)$. 
\end{proof}

\begin{corollary} \label{c:TildeCont} 
	Let $\pi \colon G \times M \to M$ be a continuous action by isometries on $(M,\rho)$ with fixed $\0$. 
	Assume that the action 
	satisfies the following (equiuniformity) condition 
	$$
	\forall \eps>0 \ \ \exists U(e) : \ \ \rho(ux,x) \leq \eps \ \ \  \text{for all} \ x \in X \ \text{and} \ u \in U(e). 
	$$
	Then the natural action $\widetilde{\pi} \colon G \times M^{\mathcal{R}} \to M^{\mathcal{R}}$ is continuous. 
%
\end{corollary}
\begin{proof} 	
For every $f \in \Lip_0(M)$ the following is true 
$$
|f(ux)-f(x)|\leq ||f||_{\Lip} \cdot \rho(ux,x).
$$ 	
Now, our assumption implies that every $f \in \Lip_0(M)$ belongs to  $\mathrm{RUC}_G(M)$ and we can apply Theorem \ref{t:positive}.1. 
\end{proof}

\begin{proposition} \label{p:positive}  
Let $G$ be a metrizable abelian topological group acting on itself by translations and $\rho$ is a bounded  
invariant topologically compatible metric on $G$. Consider the pointed metric isometric $G$-space $(M,\0)$, where $M=G \cup \{\0\}$ and $\0$ is a new $G$-fixed point with $\rho(\0,g)=\operatorname{diam}(\rho)$ for every $g \in G$. 
	Then $\d_* \colon M \hookrightarrow M^{\mathcal{R}}$ 
	is a $G$-compactification.	
%
\end{proposition}
\begin{proof}
For	every $f \in \Lip_0(M)$ and $\eps >0$ there exists a neighbourhood $U(e)$ such that the following condition holds 
$$ \forall x \in G \ \ \forall u \in U(e)  \ \ \ \  
|f(ux)-f(x)| \leq 
||f||_{\Lip} \cdot \rho(ux,x) = ||f||_{\Lip} \cdot \rho(u,e) < \eps.
$$ 
Also, $f(\0)=0=f(g \0)$.  
Hence, $f$ is a 
$\mathrm{RUC}_G$ function on $M$ in terms of Definition \ref{d:uniform}. We obtain, that $\Lip_0(M) \subseteq \mathrm{RUC}^b_G(M)$. Now, 
Theorem \ref{t:positive}.2 implies that 
$\d_*$ is a $G$-compactification. 
\end{proof}

\begin{example} \label{ex:NotCont} 
	Let $M:=(\R^2,\rho)$, where $\rho(x,y):=\min\{||x-y||,1\}$. Then we get a pointed bounded metric space with $\0=(0,0)$. Consider the compact circle group $\T$ and its isometric continuous action on $M$ by rotations around $\0$. Then 
	\begin{enumerate}
		\item 
		$\d_* \colon M \hookrightarrow M^{\mathcal{R}}$ is a $G_{\mathrm{disc}}$-compactification but not a $G$-compactification; 
		\item The orbit map $orb_{f_A} \colon \T \to \Lip_0(M)$ is not norm continuous, where $f_A(x):=\rho(A,x)$ with $A:=\Z \times \{0\} \subset \R^2$.  
	\end{enumerate}
	\begin{proof} 
		(1) 
	 The bounded function $f_A \leq 1$ belongs to $\Lip_0(M)$ 
	(because, $||f_A||_{\Lip}=1$ and $f_A(\0)=0$)
	but $f_A \notin \mathrm{RUC}_G(M)$. Indeed, take $x_n:=(n,0)$. Then $f_A(x_n)=0$ for every $n \in \N$ but for every neighborhood $U(e)$ in $\T$ there exist sufficiently big $n$ and $g_n \in U(e)$ such that $f_A(g_nx_n)=1$. This proves (1) (using Theorem \ref{t:positive}.2). 
	
	\sk 
	(2) Choose $x_n:=(n,0), y_n:=(n+\frac{1}{n},0), \ n \in \N$. Then for every neighborhood $U(e)$ there exist sufficiently big $n$ and $g_n \in U(e)$ such that $f_A(g_n^{-1}x_n)=f_A(g_n^{-1}y_n)=1$. Then $|(g_nf_A-f_A)(x_n)-(g_nf_A-f_A)(y_n)|=\frac{1}{n}=|x_n-y_n|$. Therefore, $||g_nf_A-f_A||_{\Lip} \geq 1$. 
	\end{proof} 
	\end{example}


\sk   
\subsection*{Equivariant Gromov compactifications} 

\begin{definition} \label{d:GromCom} 
	Let $(X,d)$ be a \textbf{bounded} metric space (not necessarily pointed) and $G \times X \to X$ is a continuous isometric action. Consider the following (bounded) family of distance functions 
$$
\Gamma:=\{\g_a \colon X \to \R, \ \g_a(x):=d(a,x)\}_{a \in X}.
$$
	This family $\Gamma$ is $G$-invariant ($g \g_a=\g_{ga}$) and $\Gamma \subseteq \mathrm{RUC}_G(X)$ (because $|\g_a(gx)-\g_a(x)|\leq d(g^{-1}a,a)$). 
	
	Let $Gro(X)$ be the closed unital subalgebra of $C_b(X)$ 
	generated by $\Gamma$.  
	Denote by $\g \colon X \to \widetilde{X}^{\g}$ the corresponding compactification. Remark \ref{r:equivalence} guarantees that $\g$ is a $G$-compactification.
	Following \cite{Me-opit07}, \cite[Section 2]{IbMe20},  
	we call it the \textbf{Gromov compactification} of the isometric $G$-space $X$.  
\end{definition} 

Note that $\g$ is a topological embedding because $\Gamma$ separates points and closed subsets. Indeed, for every closed subset $B \subset X$ and $x_0 \in X \setminus B$, we have $\g_{x_0}(x_0)=0$ and $\g_{x_0}(b)\geq d(x_0,B)$ for every $b \in B$. Consider the family of induced bounded pseudometrics 
\begin{equation} \label{eq:Gromov2}
	\Gamma^*:=\{\g_a^* \colon X \times X \to \R, \ \g_a^*(x,y):=|d(a,x)-d(a,y)|\}_{a \in X}. 
\end{equation} 
The corresponding weak uniformity on $X$ generates a precompact uniformity (that is the uniformity with compact completion) and its completion is just the compactification $\g \colon X \to \widetilde{X}^{\g}$. The algebra of this compactification is $Gro(X)$ as it follows by Fact \ref{f:SW}.

For some examples and applications regarding Gromov compactification we refer to \cite{IbMe20} and \cite{Pest-Smirnov}. 



\begin{remark}[space of \textit{metric types}] 
	\label{r:Gar}  	
Garling studied in \cite{Gar82} the space $T(M)$ of \textit{types} for separable metric spaces $(M,d)$. It is  a natural ``local compactification":  
$$\mathcal{T} \colon M \hookrightarrow T(M) \subset \lambda_1(M),$$ 
where $\mathcal{T}(x) \colon M \to \R$ is the distance function $\mathcal{T}(x)(y):=d(x,y)$ for every $x \in M$. Here   
 $\mathcal{T}$ is a dense topological embedding into a locally compact $\s$-compact space $T(M)$ and $\lambda_1(M) \subset \R^M$ is a topological space of all $1$-Lipschitz functions on $M$ with the pointwise topology  
 (i.e.\ the set of all real‐valued 1–Lipschitz functions on $M$).   
\end{remark}

\begin{proposition} \label{p:Garling}
Let $\pi \colon G \times M \to M$ be an isometric continuous action on a (non-pointed) separable metric space $(M,d)$. We have a natural uniquely defined action  $\pi_T \colon G \times T(M) \to T(M)$ 
on the space $T(M)$ of metric types
which extends $\pi$ and is continuous. If the metric $d$ is  bounded then the  embedding $\mathcal{T} \colon M \hookrightarrow T(M)$ is a $G$-compactification  which is equivalent to the Gromov $G$-compactification $\g \colon M \to \widetilde{M}^{\g}$. 	
\end{proposition}
\begin{proof}  
Consider the weak uniformity $\mathcal{U}$ on the set $\lambda_1(M)$ generated by the family of all projections 
$q_{x_0} \colon \lambda_1(M) \to \R, f \mapsto f(x_0)$, where $x_0 \in M$. 
The corresponding topology is just the pointwise topology on $\lambda_1(M)$.
Every restriction $q_{x_0}|_{\mathcal{T}(X)}$ on the image $\mathcal{T}(M)=M$ is the distance from $x_0$ function on $M$. Indeed, $q_{x_0}(\mathcal{T}(x))=\mathcal{T}(x)(x_0)=d(x,x_0)$. Then the natural action  
$$\pi_1 \colon G \times \lambda_1(M) \to \lambda_1(M), \ (gf)(y):=f(g^{-1}y)$$ 
is a well defined extension of $\pi$ (because $\mathcal{T}(gx)=g\mathcal{T}(x)$) with $\mathcal{U}$-uniform $g$-translations. 
The subset $\mathcal{T}(X)$ and its closure $T(X)$ are $G$-subsets and every restricted projection $q_{x_0}|_{\mathcal{T}(X)}$ is RUC in the sense of Definition \ref{d:uniform}.1  for any given $x_0 \in M$, as it follows from the following computations:  
$$
|\mathcal{T}(gx)(x_0)-\mathcal{T}(x)(x_0)|=|d(gx,x_0)-d(x,x_0)|=
|d(x,g^{-1}x_0)-d(x,x_0)|\leq d(g^{-1}x_0,x_0).
$$
Thus, the action on $\mathcal{T}(M)$ is $\mathcal{U}|_{\mathcal{T}(M)}$-equiuniform, and the induced action $\pi_T \colon G \times T(M) \to T(M)$ is jointly continuous by Fact \ref{f:G-compactFacts}.1.  

 For bounded metric $d$ the space $T(X)$ is compact and every composition $q_{x_0} \circ \mathcal{T}(x) \colon M \to \R$ is just the function $\g_{x_0}$ defined in Definition \ref{d:GromCom}. Therefore, the embedding $\mathcal{T} \colon M \hookrightarrow T(M)$ is a $G$-compactification which is equivalent to the Gromov $G$-compactification $\g \colon M \to \widetilde{M}^{\g}$.
\end{proof}

Some results of this section makes sense also for topometric $G$-spaces under natural adaptations. For simplicity we do not consider these aspects in the present paper.

\sk 
\section{Representation of dynamical systems on Lipschitz-free spaces} \label{s:Repr}  

\begin{definition} \label{d:RepresFL} 
	Let $X$ be a topological $G$-space and $M$ be a pointed metric space. A  
	{\textbf{representation}} of $(G,X)$ on the Lipschitz-free space $\mathcal{F}(M)$ is a pair $(h,\a)$ where 
	$h \colon G \to \mathrm{Is}(M)$ is a continuous homomorphism and $\a \colon X \to \mathcal{F}(M)^*$ is a weak-star continuous bounded $G$-equivariant map.  
\end{definition}

This is a particular case of Definition \ref{d:repr}. Indeed,  take into account Lemma \ref{l:emb} which asserts that  $\Is(M)$ can be treated as a topological subgroup of  $\mathrm{Is}_{\text{lin}} (\mathcal{F}(M))$.  
In contrast to Definition \ref{d:repr}, here we restrict attention to homomorphisms into the isometry group $\mathrm{Is}(M)$ of the original metric space $M$, viewed as a subgroup of the linear isometry group $\mathrm{Is}_{\text{lin}}(\mathcal{F}(M))$ via Lemma \ref{l:emb}.

\begin{question} \label{q:REPR}
	Let $\mathcal{K}$ be a certain good class of pointed metric spaces. 
	Which dynamical systems $(G,X)$ can be properly represented (in the sense of Definition \ref{d:RepresFL}) on 
	$\mathcal{F}(M)$ for some $M \in \mathcal{K}$ ? 
\end{question}

Recall that (in view of Definition \ref{d:repr}) \textit{proper representation} simply means that $\a$ is a topological embedding. 
Every (proper) representation  of a $G$-space $X$ induces a (proper)  $G$-compactification. Indeed, the weak-star closure of $\a(X)$ into the dual $V^*$ induces a $G$-compactificatioin of $X$.   

Next we are going to show that there are sufficiently many representations on Lipschitz-free spaces. 
We use a (contravariant) assignment in the spirit of the Gelfand–Naimark duality, where to a compact space we associate the Banach space \( C(K) \). 
To every $G$-factor map \( q: K_1 \to K_2 \) between compact $G$-spaces, we have the induced isometric embedding 
\[ 
q^*: C(K_2) \hookrightarrow C(K_1) 
\] 
of Banach algebras, 
which restricts to a canonical isometric embedding of (pointed) metric $G$-spaces 
\[ 
q^* \colon B_{C(K_2)} \hookrightarrow B_{C(K_1)}. 
\] 
Similarly, we have another contravariant assignment in the converse direction: 
for every non-expansive $G$-map \(\alpha \colon M_1 \to M_2 \) between pointed isometric $G$-spaces, there is an induced (see 
Fact \ref{f:facts}.5) weak-star continuous affine $G$-map 
\[
\Gamma := \mathcal{F}(\alpha)^* \colon B_{\mathcal{F}(M_2)^*} \to B_{\mathcal{F}(M_1)^*}.
\]

\begin{theorem}  \label{t:suff}  
Let $K$ be a compact $G$-space then $(G,K)$ admits a proper representation on $\mathcal{F}(M)$, where $M:=B_{C(K)}$ is the norm closed unit ball as the desired pointed metric space. 
\end{theorem}
\begin{proof} Given continuous action $G \times K \to K$,  
the standard compactness argument shows that every $f \in C(K)$ is RUC in the sense of Definition \ref{d:uniform}.1. That is, $C(K)=\mathrm{RUC}^b_G(K)$. 
This implies that the induced linear isometric action 
$$G \times C(K) \to C(K), \ (gv)(x):=v(g^{-1}x),$$ on the Banach space $(C(K),||\cdot||_{\sup})$ is separately continuous, Hence, also jointly continuous by virtue of Fact \ref{f:cont}. Therefore we have the following (restricted) isometric continuous action 
$$G \times B_{C(K)} \to B_{C(K)}.$$ 

For every $a \in K$ define $p_a \in \Lip_0(B_{C(K)})$ by 
$$p_a \colon B_{C(K)} \to \R, \ \ p_a(v):=v(a) \ \text{for every} \  v \in B_{C(K)}.$$
Then
$|p_a(v_1)-p_a(v_2)|=|(v_1-v_2)(a)|\leq 1 \cdot ||v_1-v_2||_{\sup}$ and $p_a(\0)=0$. Thus, 
$p_a \in \mathrm{Lip}_0(B_{C(K)})$ is a Lipschitz function vanishing at the origin, hence corresponds to a well-defined element of $\mathcal{F}(B_{C(K)})^*$. 
 Consider the following map 
 $$p \colon K \to (B_{\mathcal{F}(M)^*},w^*), \ p(a):=p_a.$$
  Observe that if $(a_i)$ is a net in $K$ which tends to $a$ then $v(a_i)$ tends to $v(a)$ for every continuous function $v \in M=B_{C(K)}$. 
  Hence, $(p((a_i)))$ tends to $p(a)$ in the weak-star topology  
  by Fact \ref{f:facts}.2 (because $p(K)$ is norm bounded).  
  This means that $p$ is continuous.  
   Also, $p$ is injective. Since $K$ is compact, we obtain that $p$ is a topological embedding.

We have the canonical continuous homomorphism $h \colon G \to \mathrm{Is}(B_{C(K)})=\mathrm{Is}(M)$. 
Now, observe that $p$ is $G$-equivariant. Indeed, we have 
to show that $p_{ga}=g p_a$ for every $g \in G$ and $a \in K$. 
For every $v \in M$ we have $p_{ga}(v)=v(ga)$ and also 
$$
(g p_a)(v)=p(a)(g^{-1}v)=p_a(g^{-1}v)=v(ga).
$$
One may conclude that $(h,p)$ is a proper representation (in the sense of Definition \ref{d:RepresFL}) of $(G,K)$ on the Banach space $\mathcal{F}(B_{C(K)})$.   
\end{proof}

\begin{corollary} \label{c:suff} 
	Every compact 
	$G$-flow admits a proper representation (in the sense of Definition \ref{d:RepresFL}) on some Lipschitz-free space $\mathcal{F}(M)$.
\end{corollary}


\begin{lemma} \label{l:PK-embedding}
	Let \(K\) be a compact Hausdorff space, and set
	$
	M \;:=\;\bigl(B_{C(K)},\,d),
	$
	the closed unit ball of \(C(K)\) (distinguished point \(\0\)), with the supremum metric \(d\).  Define 
	\[
	j \colon C(K)^* \;\longrightarrow\; \Lip_0\bigl(B_{C(K)},d),
	\qquad
	j(\mu)(v)\;:=\;\mu(v)\quad(v\in B_{C(K)}).
	\]
	Then:
	\begin{enumerate}
		\item \(j\) is a linear isometric embedding, i.e.\ 
		\(\|j(\mu)\|_{\Lip}=\|\mu\|\) for all \(\mu\in C(K)^*\).
		\item Its restriction to the unit ball
		\[
		j\colon \bigl(B_{C(K)^*},\,w^*\bigr)
		\longrightarrow
		\bigl(\Lip_0(B_{C(K)},d),\,w^*\bigr)
		\]
		is an affine topological embedding of the weak\(^*\)-compact ball \(B_{C(K)^*}\) into the weak\(^*\)-compact unit ball of \(\Lip_0\bigl(B_{C(K)},d\bigr)\).
	\end{enumerate}
\end{lemma}

%

\begin{proof}
	Each such function \( j(\mu) \) induces a function on \( B_{C(K)} \) which is Lipschitz with constant at most \( \|\mu\| \). 
	In fact, $j$ an isometry because  
	it is an isometry on the closed span of point‐masses. 
	The weak-star topology on \( B_{C(K)^*} \) is the topology of pointwise convergence on \( C(K) \), and in particular on the subset \( B_{C(K)} \subset C(K) \). On the other hand, by Fact~\ref{f:facts}.2, the weak-star topology on norm-bounded subsets of \( \Lip_0(B_{C(K)}, d) \) coincides with the topology of pointwise convergence on \( B_{C(K)} \). 
	Therefore, the map \( j \colon B_{C(K)^*} \to \Lip_0(B_{C(K)}, d) \) is pointwise continuous.  
	Hence, \( j \) is continuous when both domain and codomain are equipped with their respective weak-star topologies. 
	Since \( B_{C(K)^*} \) is compact in the weak-star topology, 
	it follows that the injective map \( j \) is a topological embedding onto its image. 
\end{proof}

 \begin{remark} \label{r:P(K)}
	In the setting of Theorem \ref{t:suff} we have an additional advantage regarding topological embedding  
	$p \colon K \to (B_{\mathcal{F}(M)^*},w^*)$, where $M=B_{C(K)}$.  
	As it follows from Example~\ref{e:measures} and Lemma   \ref{l:PK-embedding} 
	(applied to the usual setting of Lipschitz-free spaces), 		
	not only the space \( K \) embeds topologically into the dual space $\mathcal{F}(M)^*$ via the natural identification with point evaluations. More importantly, it is true also for the weak-star closure of the affine hull of \( K \) which is just the space \( P(K) \) of all regular Borel probability measures on \( K \).  
	This observation is particularly might be useful in the context of (metric) amenability of group actions.

%
%
\end{remark}






 \begin{definition} \label{d:affCOMP} 
 	Let $(h,\a)$ be a continuous representation of a $G$-space $X$ on a Banach space $V$. 
 	\begin{enumerate}
 		\item We say that this representation is \emph{amenable} if the induced affine $G$-compactification (Definition \ref{d:affineComp}) of $X$, defined by $Q_{\a}:=\overline{co}^{w^*}(\a(X))$,  contains a $G$-fixed point. 
 		\item 	 We say that $f \in V^*$ is \emph{amenable} if $Q_f:=\overline{co}^{w^*}(Gf)$ contains a $G$-fixed point. Or, equivalently, if the tautological representation of the $G$-space $X:=Gf$ on $V$ is amenable.    
 	\end{enumerate} 
 \end{definition}
 
 For every abstract compact $G$-flow $K$ we have a canonical proper representation of $(G,K)$ on $C(K)$. Then the  amenability of this particular representation in the sense of Definition \ref{d:affCOMP} is equivalent to the existence of a fixed point in the space $P(K)$ of all probability measures on $K$. That is, this agrees with one of the most traditional definitions of amenability for actions. Recall that $\overline{co}^{w^*}(K)=P(K)$, where $K$ is naturally embedded into $C(K)^*$ and any affine $G$-compactification is an affine $G$-factor of $P(K)$.   
 
 \sk 
To every $f \in {\mathcal{F}(M)^*}=\Lip_0(M)$ 
we have a canonically defined compact 
 dynamical $G$-system 
$$
K_f:=cl_{w^*}(Gf)=\overline{{Gf}}^{w^*}.
$$ 
If $M$ is separable, then $\mathcal{F}(M)$ is separable and every $K_f$ is metrizable. 
Dynamical complexity of such natural $G$-flows leads to a complexity hierarchy for Lipschitz functions on $M$. 
It seems to be an attractive task to clarify when $K_f$ is dynamically small. 

The following two questions are closely related.  
For the definitions of the algebras   
$$\mathrm{WAP}(G) \subseteq  \mathrm{Asp}(G) \subseteq \mathrm{Tame}(G),$$   
classes of $G$-flows:  
\{WAP (weakly almost period)\} $\subseteq$ 
\{HNS (hereditarily non-sensitive)\} $\subseteq$ \{tame\},    
and their roles in Banach representations theory, see, for example \cite{GM-fixp12,GM-AffComp13} and \cite{FGN24}. 

\begin{question} \label{q:ind}
	For \(f\in\Lip_0(M)\), consider \(K_f:=\overline{Gf}^{w^*}\).
	\begin{itemize}
		\item[(a)] When is \((G,K_f)\) WAP, HNS, tame?  
		Equivalently (if \(M\) separable) when does \((G,K_f)\) admit a proper
		representation on a reflexive/Asplund/Rosenthal space?
		\item[(b)] When does the affine compactification
		\(Q_f=\overline{co}^{w^*}(Gf)\) contain a fixed point?
	\end{itemize}
\end{question}

%
%
 
 
  In Theorem \ref{t:WAP} below we have a particular case with WAP $K_f$.  
Note that if $K_f$ is a WAP dynamical system, then $Q_f$ is a WAP affine (weak-star compact) $G$-flow, and $f$ is amenable. 
Moreover, $f$ (even any $\varphi \in Q_f$) is amenable already under a weaker assumption when the $G$-flow $K_f$ is only HNS 
as it follows by the following 
result. 


\begin{fact} \label{f:measure} \cite[Proposition 2.2 and Theorem 1.5]{GM-fixp12} 
Every hereditarily non-sensitive (HNS) (e.g., WAP) compact $G$-system $K$ admits an invariant probability
measure. Hence, every affine $G$-compactification of $K$ contains a $G$-fixed point and every Banach representation of $(G,K)$ is amenable. 	
\end{fact}

 Note that \cite[Theorem 1.5]{GM-fixp12} simultaneously generalizes classical Ryll-Nardzewski's fixed point theorem and also Fact \ref{f:measure}. 
 If $K_f$ is norm-separable, then any $\varphi \in Q_f$  (e.g., $f$) is amenable by a known folklore fixed-point theorem. A direct prove can be found in a work of Glasner \cite[Theorem 1.2]{Gl-YomDin}.  
Another proof can be derived from 
Fact \ref{f:measure} 
because norm separable $K_f$ is (weak-star,norm) fragmented and in this case the $G$-flow $K_f$ is HNS.  

\sk 
For every $f \in {\mathcal{F}(M)^*}$ 
and $v \in M$,  we may consider the corresponding \textbf{matrix coefficient} 
$$mat_{f,v} \colon G \to \R, \ g \mapsto f(g^{-1}v),$$ 
which is a bounded right uniformly continuous function on $G$; see \cite{Me-FRepres03}.  
  \begin{question} \label{q:WAP}   
  	When $mat_{f,v}$ belongs to a dynamically interesting class of functions? 
  	For instance, when $mat_{f,v}$ belongs to $\mathrm{WAP}(G)$, $\mathrm{Asp}(G)$, $\mathrm{Tame}(G)$ ?
  \end{question}


\begin{remark} \label{r:functorial} 
	One source of general examples is the functorial assignment from Theorem \ref{t:suff}.  
	Suppose that \( K \) is a compact \( G \)-flow containing a dense orbit \( Gx_0 \). 
	Then, in the setting of Theorem \ref{t:suff} (where \( M := B_{C(K)} \)), we obtain a proper representation of the \( G \)-flow \( K \) on the corresponding Lipschitz-free space \( \mathcal{F}(M) \).  
	In this case, \( x_0 \) can be treated as an element \( f \in \Lip_0(M) \), and \( K = K_f \). 
\end{remark}

%

\begin{proposition} \label{p:MeGl} 
	There exists an isometric action of the rank $2$ free discrete group $G=F_2$ on a separable metric space $M$ such that a certain compact $G$-subspace $K$ of $(\mathcal{F}(M)^*, w^*)$ is a tame minimal $G$-flow but the tautological representation of $(G,K)$ on $\mathcal{F}(M)$ is not amenable in the sense of Definition \ref{d:affCOMP}.  	
\end{proposition}
\begin{proof} According to \cite[Example 3.1]{GM-fixp12}, 
	there exists a tame minimal compact metric $G$-system $K$, with $G=F_2$, such that $P(K)$ does not have a fixed point. 	
	Now, apply the method of Remarks \ref{r:P(K)} (based on  Lemma \ref{l:PK-embedding}) 
	and Remark \ref{r:functorial} to the pointed metric space $M:=B_{C(K)}$ with the induced continuous isometric action of $G$ on $M$. 	
\end{proof}
                                      
\sk 
\section{Equivariant metric (horo) compactifications} 
\label{s:horo} 

We consider the so-called \textit{metric compactifications} (\textit{horocompactifications}) $\mu \colon M \to \widetilde{M}$ which is well known in metric geometry. 
There are several (different) definitions in the literature. One of the main versions of this concept was introduced by M. Gromov \cite{Gro81}. Relevant information about metric 
(horo) compactifications can be found, for example, in \cite{Gut,GutDiss,Duchesne,Daniil}.  
Following the approach of these papers briefly recall the definition.  
Let $(M,d,\0)$ be a pointed metric space. Consider the function  
$$\mu \colon M \to \R^M, \quad a \mapsto \mu_a, \quad \mu_a(x) := d(a,x) - d(a,\0).
$$ 
Here $\R^M$ carries the pointwise (product) topology. 
It is well known and easy to see that $\mu$ is always continuous and injective. The pointwise closure $\widetilde{M}:=cl(\mu(M))$ in $\R^X$ of the image is compact. Thus, we have an induced  compactification map which also will be denoted by $\mu$. This compactification map $\mu \colon M \to \widetilde{M}$ is the \textit{metric compactification} (\textit{horocompactification}, \textit{Busemann compactification}) of $(M,d,\0)$.  
The remainder $\partial (\widetilde{M}):=\widetilde{M} \setminus M$ is called the \textit{horofunction boundary}.  

Metric compactification is independent (up to  homeomorphisms) of the choice of base point.

\begin{remark} \label{r:NotEmb} 
In general, $\mu$ is not a topological embedding. See \cite[p. 25]{GutDiss} or \cite{Daniil} with $M:=(\ell_1,\0)$, 
where $\0$ denotes the zero vector in the Banach space 
$\ell_1$. Indeed, let $v_n$ be the sequence 
$v_n: = (0, ..., 0, n, 0, ...)$, with $n$ in the $n$-th coordinate. 
Then observe that $\lim_{n \to \infty} \mu(v_n)=\mu_{\0}$ but $\lim ||v_n -\0||=\infty$. 

M.I. Garrido (one of the authors of \cite{Daniil})  
informed us 
that $\mu$ need not be an embedding also for \textbf{bounded} metric spaces. Namely, for the metric subspace $M:=\{\0\} \cup \{e_n: n \in \N\}$ of $\ell_1$. 
\end{remark}

 It is well-known (see, for example, \cite{GutDiss,Daniil}) that $\mu$ is a topological embedding for every complete geodesic and \textit{proper} (meaning that all closed balls are compact) metric space $(M,d)$.  
 
 Fioravanti  proved in \cite[Proposition 4.21]{Fioravanti} that when $M$ is a complete, locally convex median space with  compact intervals, its Roller compactification  coincides up to homeomorphism with the metric (Busemann) compactification.

\begin{definition}
	Let us say that a point $x_0$ in $(M,d)$ is \textit{equidistant} 
	(or, \textit{$c_0$-equidistant}) 
	if $d(x,x_0)=c_0 >0$ is constant for 
	every $x \in M \setminus \{x_0\}$.
\end{definition} 

 Clearly, if there exists a $c_0$-equidistant point, then 
 $\operatorname{diam}(M,d) \leq 2 c_0$. 
  Conversely, for bounded metrics, one may adjoin a new point $\0$ which is equidistant as we observed in Remark \ref{r:IsoFix} 
  (take, for example, $c_0=\operatorname{diam}(M,d)$). 
  This is useful in view of actions because in this way any isometric $G$-action on a bounded metric $G$-space $M$ can be naturally embedded into an isometric $G$-action on the pointed space $M \cup \{\0\}$ fixing the new isolated point. In this case, assertion (2) of Theorem \ref{t:hori}, in fact, speaks about the ``original" non-pointed metric space $M$ and its Gromov compactification (see \cite[Prop. 2.7]{IbMe20}).



 \sk 
\begin{theorem} \label{t:hori} 
	Let $(M,d,\0)$ be a pointed isometric $G$-space, and let $h \colon G \to \Is(M)$ denote the corresponding homomorphism. 
	 Define the map 
	$$
	\mu \colon M \to (\mathcal{F}(M)^*, w^*), \  \ \mu(a)=\mu_a,
	$$	 
		$$\mu_a(x):=d(a,x)-d(a,\0).$$  
	\begin{enumerate} 
		\item 
		\begin{enumerate}
			\item The pair $(h,\mu)$ is a continuous injective representation of the $G$-space $M$ on 
			$\mathcal{F}(M)$, with 
			 $||\mu(a)||_{\mathrm{Lip}}=1$ for every $a \in M$. 
			\item The induced continuous (injective) $G$-compactification   
			$$\mu \colon M \to \widetilde{M}:= \overline{\mu(M)}^ {w^*} \subset B_{\mathcal{F}(M)^*}$$
			is equivalent to the {metric (horo)compactification} of $(M,d,\0)$. 
		\end{enumerate}
		
		 \item Let $\0$ be equidistant in $M$ with $c_0:=d(x,\0)$ for every $x \in X:= M \setminus \{\0\}$.  
		 Then 
		 \begin{enumerate}
		 	\item the restriction map   
		 	$$\mu|_X \colon X \to  \overline{\mu(X)}^{{w^*}}$$ 
		 	is a proper (topological embedding) compactification  and is equivalent to the Gromov compactification (Definition  \ref{d:GromCom}) of $(X,d)$.  
		 	\item If $\operatorname{diam} (M \setminus \{\0\}) < 2c_0$, then 
		 	$\mu \colon M \to \widetilde{M}$ is 
		 	a topological embedding.  
		 \end{enumerate} 
	\end{enumerate}    
\end{theorem}
\begin{proof}
	We repeatedly use the equality $\mathcal{F}(M)^*=\mathrm{Lip}_0 (M)$ (Fact \ref{f:facts}.1). 
	
	 (1) It is straightforward to see that $\mu_a \in \mathrm{Lip}_0 (M)$,  $||\mu_a||_{\mathrm{Lip}}=1$ for every $a\in M$ and	\textbf{$\mu$ is} \textbf{injective}.  
	
	
%
	
	 	\sk  
	 \textbf{$\mu$ is weak-star \textbf{continuous}.} 
	 For every $x \in M$,  
	 define the following function 
	 $$\varphi_x \colon M \to \R, \  \varphi_x(a):=d(a,x)-d(a,\0)=\mu_a(x).$$
	 Then $\varphi_x$ is bounded because $|\varphi_x(a)| \leq d(\0,x)$ and continuous (being $2$-Lipschitz) by  
	 $$
	 |\varphi_x(a_1)-\varphi_x(a_2)| \leq 2 d(a_1,a_2).
	 $$ 
	  Every $\mu_a$ can be identified with the following element of $\R^M$ defined as follows: 
	 $$(\mu_a(x))_{x \in M}= (d(a,x)-d(a,\0))_{x \in M}=(\varphi_x(a))_{x \in M} \in \R^M.$$
	 Since $\varphi_x(a)=\mu_a(x)$, the function 
	 $\mu \colon M \to \R^M, \ a \mapsto \mu_a$
	 is the diagonal product of the following family of functions $\varPhi:=\{\varphi_x: \ x \in M\}$. 
	 Denote by $\tau_w$ the corresponding pointwise (weak) topology on $\mu(M)$ which coincides with the topology of the corresponding precompact uniformity $\mu_\varPhi$ on $M$.

	 As we have already established, 
	  $\mu(M) \subset B_{\mathcal{F}(M)^*} \subset \mathrm{Lip}_0 (M)$ holds.  
	  Since $\mu(M)$ is norm bounded in $\mathrm{Lip}_0 (M)$, the weak-star topology inherits on $\mu(M)$ the pointwise topology (Fact \ref{f:facts}.2).  That is exactly the subspace topology $\tau_w$ of the product $\R^M$.  
	 Every $\varphi_x \colon M \to \R$ is continuous. Hence, $\tau_w \subseteq top(d)$. This implies that the injection $\mu \colon M \to B_{\mathcal{F}(M)^*}$ is weak-star  continuous. 
	 
	 Furthermore, by Fact \ref{f:SW} the metric compactification 
	 $m \colon M \to \widetilde{M}$ is the completion of the precompact uniformity $\mu_\varPhi$ on $M$ generated by the family of (bounded $2$-Lipschitz) functions 
	\begin{equation} \label{eq:1} 
		\varPhi:=\{\varphi_x \colon M \to \R, \ \varphi_x (a)=d(a,x)-d(a,0))\}_{x \in M}.
	\end{equation}
	 	
	In other words, the topology of $\mu(M)$ inherited from $\widetilde{M}$ is the weak topology  
	generated by the family $\varPhi$. 
	
	 	\sk  
	 	\textbf{$\mu$ is $G$-equivariant.} 
	 	That is, $\mu(g a)=g \mu(a)$ for every $g \in G$. 
	 	Indeed, taking into account the description of the dual action (see Definition \ref{d:repr}), for every $x \in M$ we obtain  
	 \begin{align*}
	 	\mu(ga)(x)=\mu_{ga}(x)=d(ga,x)-d(ga,\0)	=d(a,g^{-1}x)-d(a,g^{-1}\0) \\
	 	=d(a,g^{-1}x)-d(a,\0)=\mu_a(g^{-1}x)=(g\mu_a)(x)
	 \end{align*}
	 
	 	\sk 
	 	
	  (2a)   
	   We have to show that $\mu|_X$ is a 	
	  \textbf{topological embedding} of $X:=M \setminus \{0\}$  into the weak-star compact space $(B_{\mathcal{F}(M)^*},w^*)$, assuming that $c_0:=d(x,\0)$ for every $x \in X:= M \setminus \{\0\}$. 
	   Obviously, $\varphi_x (a)=d(a,x)-c_0$ for every $a \in X$ and every $x \in X$.
	  By Fact \ref{f:SW} and a discussion above before (\ref{eq:1}), it is enough to show that the following family of functions 
	  $$
	  	\Gamma_0:=\{\varphi_x \colon X \to \R, \ \varphi_x (a)=d(a,x)-c_0\}_{x \in X} 
	  $$
	   separates points and closed subsets in $X$. Given any closed set $B \subset X$ and point $x_0 \in X \setminus B$, we note that $\varphi_{x_0}(x_0)=-c_0$, while for each 
	   $b \in B$ we have   
	  $\varphi_{x_0}(b)\geq c_1 - c_0$, where $c_1:=d(x_0,B)>0$. Hence, $\varphi_{x_0}(x_0) \notin cl(\varphi_{x_0}(B))$. 
	  
	  Clearly, the family $\varPhi_X:=\{\varphi_x|_{X}: \ x \in X\}$, with $X:=M \setminus \{\0\}$, generates the same unital subalgebra $Gro(X)$ of $C_b(X)$ as the family 
	  $$\Gamma:=\{\g_a \colon X \to \R, \ \g_a(x):=d(a,x)\}_{a \in X}$$
	 from Definition \ref{d:GromCom}. This implies that $\mu|_X \colon X \to  \overline{\mu(X)}^{{w^*}}$ is equivalent to the Gromov compactification 
	 of $(X,d)$.
	 
	 \sk  	 
	 (2b) 
	 Now, assume, in addition, that $\operatorname{diam} (X) < 2c_0$. By (2a), $\mu|_X \colon X \to  \overline{\mu(X)}^{{w^*}}$ is a topological embedding. 
	 The point $\0$, being equidistant, is an isolated point in $M=X \cup \{\0\}$. Thus, it is enough to show that some $\varphi_y$ separates $\0$ and $X$. Take any $y \in X$. Then 
	 $$\varphi_y(\0)=d(\0,y)-d(\0,\0)=d(\0,y)=c_0$$ and 
	 $$\varphi_y(x)=d(x,y)-d(x,\0)=d(x,y)-c_0 \leq \operatorname{diam}(X)-c_0.$$  
	 Since $\operatorname{diam} (X)-c_0 < c_0$, we obtain   
	 $$
	 \varphi_y(x) \leq \operatorname{diam}(X)-c_0 < \varphi_{y}(\0),  
	 $$
	 for every $x \in X$. Hence, $\varphi_{y}(\0) \notin  cl(\varphi_{y}(X))$.  
\end{proof}

The continuity of the induced $G$-action on $\widetilde{M}$ was verified in \cite[Lemma 2.5]{Duchesne}. This fact also follows directly from Theorem \ref{t:hori} and Lemma \ref{l:DualAction}. 
Assertion (2) of Theorem \ref{t:hori}
shows that the Gromov compactification yields interesting geometric examples of representations (in the sense of Definition \ref{d:RepresFL}) on Lipschitz-free spaces.

\begin{remark}  \label{r:link} 
Metric compactification $\widetilde{M}$ is a natural factor of 
the Lipschitz realcompactification $M^{\mathcal{R}}$ (see Remark \ref{r:APS}). Indeed, recall that $\d_* \colon M \hookrightarrow M^{\mathcal{R}}$ 
	arises as the completion of the weak uniformity $\mathcal{U}_*$ which comes on $M$ from the family of functions $\Lip_0(M)$.  
	Since the family $\varPhi$ from (\ref{eq:1}) is contained in $\Lip_0(M)$, there exists a continuous onto map $$q \colon M^{\mathcal{R}} \to \widetilde{M}.$$
	Now, if $M$ is a $G$-space under an isometric action then $\varPhi$ is $G$-invariant and $q$ is equivariant. 

If, in addition, $\Lip_0(M) \subseteq \mathrm{RUC}_G(M)$ (as in Theorem \ref{t:positive}.2), then $M^{\mathcal{R}}$ is a $G$-space, and when $M$ is bounded, the map $q$ becomes a $G$-compactification factor.

\end{remark}


\sk  
We say that a map $F \colon A \times B \to \R$ has the \textit{Double Limit Property} (in short: DLP) if for every pair of sequences 
$(a_n)_{n \in \N}, (b_m)_{m \in \N}$ in $A$ and $B$ respectively, 
$$
\lim_n \lim_m F(a_n,b_m) \ = \ \lim_m \lim_n F(a_n,b_m)
$$ 
whenever both of these limits exist. 
In particular, for the map $d \colon M \times M \to \R$, this gives a well-known definition (see, \cite{Gar82}) of the \textit{stable metric} $d$ which is a natural generalization of stable norms. Let $G \times X \to X$ be a group action. We say that $f \colon X \to \R$ has the DLP if the induced map 
 $$
 w_f \colon Gf \times X \to \R, \ (gf,x) \mapsto f(g^{-1}x)
 $$ 
 has the DLP.  
A continuous bounded function $f \colon X \to \R$ defined on a $G$-space $X$ is WAP if $f$ has the DLP. This DLP-criterion of WAP (in fact, relative weak compactness of $Gf$ in $C_b(X)$) is well known and goes back to Grothendieck \cite[Ch.,X, Thm X]{BJM}. See, for example, 
\cite{BJM} and \cite[Fact 2.4 and Theorem 8.5]{Me-FRepres03}.

\begin{theorem} \label{t:WAP} 
Let $(M,d)$ be a bounded pointed metric space with an isometric continuous $G$-action. Suppose that $d$ is a stable metric. Then 
\begin{enumerate} 
	\item The metric $G$-compactification $\widetilde{M}$ is a WAP $G$-flow.   
	\item The $G$-flows $\widetilde{M}$ and $K_{\mu_a}$ admit  proper representations on reflexive Banach spaces for every separable $M$ and $a \in M$. Every functional $\mu_a$ is amenable. 
		\item $mat_{\mu_a,v} \in \WAP(G)$ for every $a,v \in M$. 
\end{enumerate} 	
\end{theorem} 
\begin{proof} 
	(1) 
	Recall (see (\ref{eq:1})) that the following family of bounded Lipschitz functions 
	 \begin{equation} \label{eq:GromovNewPar} 
		\Gamma_0:=\{\varphi_z \colon M \to \R, \ \varphi_z (x)=d(x,z)-d(x,\0))\}_{z \in M} 
	\end{equation} 
	generates the metric compactification $\mu \colon M \to \widetilde{M} \subset \R^M$. Observe that $\Gamma_0$ is $G$-invariant because $g \varphi_z=\varphi_{gz}.$ The corresponding algebra of the compactification $\mu$ is the smallest Banach subalgebra of $\RUC_G(M)$ containing $\Gamma_0$ and constants as it follows by Fact \ref{f:SW}. 
	
	Always, $\WAP(X)$ is a $G$-invariant norm closed subalgebra of $C_b(X)$ for every $G$-space $X$. Thus, it is enough to show that every $\varphi_z$ belongs to $\WAP(M)$. This follows by above mentioned DLP-criterion of WAP.   
In order to verify that $\varphi_z \in \WAP(M)$, we have to show that the induced map 
	$$w \colon G\varphi_z \times M \to \R, \ \ w(g\varphi_z,x):=\varphi_z(g^{-1}x)$$ 
	has the DLP. Observe that 
	$$
	\varphi_z(g_n^{-1}x_m)=d(g_n^{-1}x_m,z)-d(g_n^{-1}x_m,\0)=
	d(x_m,g_nz)-d(x_m,g_n\0)=d(x_m,g_nz)-d(x_m,\0).
	$$
	Since the double sequence $d(g_n^{-1}x_m,z)-d(g_n^{-1}x_m,\0)$ is bounded, one may suppose, up to passing to subsequences (see, \cite[Lemma 3.3]{Me-Geo04}) that there exist the corresponding double limits. Moreover, since $d$ is bounded, one may suppose, in addition, that there exists $\lim_{m \in \N} d(x_m, \0)=r \in \R.$ 
	Then  
	$$
	\lim_n \lim_m (d(x_m,g_nz)-d(x_m,\0))=\lim_n \lim_m d(x_m,g_nz)-r
	$$ 
	$$\lim_m \lim_n (d(x_m,g_nz)-d(x_m,\0))=\lim_m \lim_n d(x_m,g_nz)-r
	$$ 
	Finally, use the DLP (stability) of $d$.  

	\sk 
	(2) 
	By (1), $(G,\widetilde{M})$ is a WAP $G$-flow. If $M$ is separable, then also $\mathcal{F}(M)$ is separable. Thus, the compact space $(B_{\mathcal{F}(M)^*},w^*)$ is metrizable. Therefore, $\widetilde{M}$ is also a metrizable compact $G$-flow. Now,  $\widetilde{M}$ (being a metrizable WAP $G$-flow)  admits a proper representation on a reflexive Banach space by \cite{Me-FRepres03}.  
	The same is true for $K_{\mu_a}$ because it is a $G$-subflow of $\widetilde{M}$. Finally, $\mu_a$ is amenable by Fact \ref{f:measure}.   
	
	(3) As in (1), we use the DLP of $d$ (and the DLP criterion of Grothendieck) taking into account the following  equality: 
	$$mat_{\mu_a,v}(g_nh_m)=d(a,(g_nh_m)^{-1}v)-d(a,\0)=
	d(h_ma,g_n^{-1}v)-d(a,0).$$  
\end{proof}

 A Banach space $(V,||\cdot||)$ is said to be \textit{stable} (Krivine and Maurey \cite{KrMa}, and \cite{Gar82}) if the natural norm metric is stable. It is well known that all $L_p(\mu)$ Banach spaces 
are stable for every $1 \leq p <\infty$.

\begin{corollary} \label{c:StableBanach} 
For every stable Banach space $(V,||\cdot||)$ (e.g. $V:=L_p(\mu)$) and the natural isometric action of the topological group $G:=\Is_{\text{lin}}(V)$ on the closed unit ball $B_V$ (with the origin $0_V$ as the distinguished point) the corresponding metric $G$-compactification $\widetilde{B_V}$ is a WAP $G$-flow. Every metric functional $\mu_a$ is amenable for each $a \in B_V$. 
\end{corollary}




\sk  
Discussion above regarding hereditary non-sensitive $G$-flows property after Fact \ref{f:measure}, and also Question \ref{q:MetrAmen} below, justify the following query. 

\begin{question} \label{q:tame} 
	For which 
	metric $G$-spaces $M$ the metric $G$-compactification $\widetilde{M}$ is a hereditarily non-sensitive (or, at least, tame) $G$-flow ? For which Banach spaces $V$ this happens for the pointed metric space $M:=B_V$ ? 
\end{question}



\begin{remark} \label{r:DouchaK} 
Question \ref{q:tame} makes sense in a more general setting for all isometric $G$-actions, where $M$ is not necessarily a pointed space. Note that quite often induced isometric actions on $\widetilde{M}$   have a fixed point in the horofunction boundary $\partial (\widetilde{M})$. Such fixed points often arise when the group action admits contracting or quasi-invariant directions at infinity, as seen in CAT(0) spaces or hyperbolic-type geometries.
See, for example \cite{Karlsson}. We thank M. Doucha, who advised us this work of A. Karlsson. 	
\end{remark}

\begin{remark} \label{r:amenability} 
Recall that one of the most common definitions of amenability for general topological groups $G$ is the existence of a fixed point in every affine compact $G$-flow (Day's fixed point theorem).  
See, for example, \cite[Theorem III.3.1]{Gl-book} and \cite[Theorem G.1.7]{BHV}). 
It is equivalent that every $G$-flow $K$ admits an invariant probability measure $m \in P(K)$. 
 That is, the affine $G$-flow $P(K)$ has a $G$-fixed point.   
 One more equivalent conditions is: the maximal $G$-compactification $\beta_G G$ of $G$ (with respect to the left action) is amenable (Definition \ref{d:affCOMP}).  
 
The Lipschitz-free setting and metric geometry suggest examining a weaker kind of amenability with a certain metric flavor (besides Question \ref{q:ind}.b see Question \ref{q:MetrAmen} below). 
\end{remark}

%

\begin{remark} \label{r:affine} (\emph{Affine  horocompactification}) 

It is natural to extend the horocompactification by considering not just individual metric functionals, but their convex combinations, thereby passing to a more flexible affine-geometric framework.
	
Let again $(M,d,\0)$ be a pointed isometric $G$-space and  
$$\mu \colon M \to \widetilde{M}:= \overline{\mu(M)}^ {w^*} \subset B_{\mathcal{F}(M)^*}$$ 
is its metric $G$-horocompactification. 
Consider the weak-star closed affine envelope 
$$\widetilde{M}^{\mathrm{aff}}:
= \overline{co}^{w^*}(\widetilde{M}) \subset \Lip_0(M),$$ 
 Then $\widetilde{M}^{\mathrm{aff}}$ is a weak-star compact convex subset of the unit ball $B_{\mathcal{F}(M)^*}$. The dual action of $G$ on $\widetilde{M}^{\mathrm{aff}}$ is continuous by Lemma \ref{l:DualAction}. That is, we get an affine $G$-compactification of $M$. 
We call to  $\widetilde{M}^{\mathrm{aff}}$ the \emph{affine horocompactification} and to 
$\widetilde{M}^{\mathrm{aff}} \setminus {M}$ the \emph{affine horofunction boundary}.   
\end{remark} 


\begin{question} \label{q:MetrAmen} 
	What is the role of 
	Lipschitz functions $f \in \widetilde{M}^{\mathrm{aff}} \setminus \widetilde{M}$ and of extreme points of $\widetilde{M}^{\mathrm{aff}}$ in the theory of horocompactifications? 
	In particular, under which conditions does there exist a $G$-fixed point in $\widetilde{M}^{\mathrm{aff}} \setminus M$? 
	
	In the latter case, it is natural to say that the original $G$-action on $M$ is \textbf{metrically amenable}. 
	
Which $f \in \widetilde{M}^{\mathrm{aff}} \setminus M$ are amenable ? 
\end{question}

One sufficient (by Theorem \ref{t:WAP}) condition is the case where the $G$-flow $\widetilde{M}$ is HNS (e.g., $(M,d)$ is  stable and bounded); see Question \ref{q:tame}. 

\begin{remark} \label{r:nonpointed}
	One may also develop a \emph{non-pointed} theory, replacing 
	\(\Lip_0(M)\) by \(\Lip(M)/\{\text{consts.}\}\).  This avoids choosing a basepoint and leads naturally to the study of metric quotients of the free space—see Weaver \cite{Weaver}.  Another approach, by Cúth–Doucha \cite{CD23}, works even when the distinguished point is not \(G\)\nobreakdash-fixed, embedding the resulting affine action into a continuous linear one.
\end{remark}

\sk 
\bibliographystyle{plain}

\end{document}